\documentclass[1p]{elsarticle}

\usepackage{lineno,hyperref}
\modulolinenumbers[5]

\journal{Computer-Aided Design}

\usepackage{amsfonts}
\usepackage{amssymb}
\usepackage{amsbsy}
\usepackage{amsmath}
\usepackage{amsthm}

\usepackage{exscale}
\usepackage{fontenc}
\usepackage{ifthen}
\usepackage{latexsym}
\usepackage{syntonly}
\usepackage{caption}
\usepackage{float}
\usepackage{array}
\usepackage{subfigure}
\usepackage{fancyhdr}
\usepackage{enumitem}
\usepackage{algpseudocode}
\usepackage{algorithm}
\usepackage{tikz}
\usetikzlibrary{matrix}
\usepackage{tikz-cd}
\usetikzlibrary{trees}
\usetikzlibrary{arrows}
\usepackage{scalerel,stackengine}

\usepackage{pifont}

\usepackage{multirow}

\def\eps{\varepsilon}
\def\bu{\mathbf u}

\renewcommand{\mathbf}[1]{\mbox{\boldmath{$#1$}}}

\newcommand{\xComplexity}{\mathcal{C}}

\newcommand{\xMat}[1]{\mathbf{#1}}  
\newcommand{\vect}[1]{\mathbf{#1}}                                    

\newcommand{\mbf}{\mathbf}

\newcommand{\parcial}[2]{\frac{\partial#1}{\partial#2}}

\newcommand{\xdeviationLinear}{\xdistortionLinear}
\newcommand{\xdistortionLinear}{\eta}

\newcommand{\xQuality}[1]{q_{#1}}

\newcommand{\xnodeBasic}{\vect{z}}
\newcommand{\xnode}[1]{\xnodeBasic_{#1}}

\newcommand{\xElement}{E}

\newcommand{\xparametrization}{\boldsymbol{\varphi}}
\newcommand{\xparametrizationDiscrete}{\xparametrization_h}
\newcommand{\xparametrizationP}{\xparametrization_{\xOrderHO}}
\newcommand{\xcurve}{\boldsymbol{\gamma}}

\newcommand{\xParamPlane}{\Pi_{\xSurface}}
\newcommand{\xSurface}{\Sigma}

\newcommand{\xreal}[1]{\mathbb{R}^{#1}}

\newcommand{\xTransposed}{^t}
\newcommand{\xParamPoint}[1]{\vect{x}_{#1}}
\newcommand{\xParamCoordFirst}[1]{x_{#1}}
\newcommand{\xParamCoordSecond}[1]{y_{#1}}
\newcommand{\xParamCoord}[1]{(\xParamCoordFirst{#1},\xParamCoordSecond{#1})\xTransposed}

\newcommand{\xMesh}{\mathrm{M}}

\newcommand{\xidealElem}[1]{\xElement^{I}_{#1}}
\newcommand{\xIdealElem}[1]{\xidealElem{#1}}

\newcommand{\xPhysicalCoord}[1]{\xnode{#1}}
\newcommand{\xPhysicalCoordFirst}[1]{x_{#1}}
\newcommand{\xPhysicalCoordSecond}[1]{y_{#1}}

\newcommand{\xrepresentationIP}{\boldsymbol{\phi}}

\newcommand{\xdeviationHO}{\eta_{{}}}
\newcommand{\xDistortionHO}{\xdeviationHO}
\newcommand{\xDistortionHOind}[1]{\eta_{{#1}}}

\newcommand{\xOrderHO}{p}
\newcommand{\numberNodesHO}[1]{{{n}_{#1}}}

\newcommand{\xDifferential}[1]{\xMat{D}#1}

\newcommand{\xNumElements}{N_{E}}
\newcommand{\xNumNodes}{N}

\newcommand{\xNumNodesLoop}{N_l}
\newcommand{\xNumElementsNeigh}[1]{N_{E}^{#1}}

\newcommand{\xNumNodesSurface}{N_{N_s}}

\newcommand{\xbnabla}{\mbf{\nabla}}

\newcommand{\xgrad}{\xbnabla}

\newcommand{\xExtraConstant}{}

\newcommand{\xDim}{d}

\newcommand{\xelind}{e}

\DeclareMathOperator*{\argmin}{argmin}

\newcommand{\xPseudoNormal}{\vect{n}_{\xnode{}}}

\newcommand{\tetgen}{{TetGen}}
\newcommand{\alya}{{Alya}}

\newcommand{\xMetric}{\mathcal{M}}
\newcommand{\xMetricParam}[1]{\xMetric_{#1} }
\newcommand{\xMetricFisrtFF}{\xMetricParam{1}}
\newcommand{\xMetricTangent}{\xMetricParam{T}}
\newcommand{\xMetricCurvature}{\xMetricParam{\xHessianNice}}
\newcommand{\xMetricCurvatureComplexity}{\xMetricParam{C}}

\newcommand{\xHessianNice}{\mathcal{H}}

\newcommand{\zh}{\mathrm{z}_h}
\newcommand{\zp}{\mathrm{z}_{\xOrderHO}}

\newcommand{\xlengthMT}{l_{\xMetricTangent}}
\newcommand{\xlengthMC}{l_{\xMetricCurvatureComplexity}}
\newcommand{\listElemsToRefine}{elemsToRefine}

\newcommand{\xNewNormSpace}[2]{ {{\lVert} #1 {\rVert}}_{#2}}

\newcommand{\xNewSum}{\sum} 

\newcommand{\xWindMesh}{{W}ind{M}esh}

\graphicspath{{./figures/}}

\usepackage{atbegshi,picture}

\AtBeginShipout{\AtBeginShipoutUpperLeft{%
		\put(\dimexpr\paperwidth-1cm\relax,-1.5cm){\makebox[0pt][r]{
				\framebox{
				}
				}
}}

\begin{document}

\begin{frontmatter}

\title{A hybrid meshing framework adapted to the topography to simulate Atmospheric Boundary Layer flows}

\date{January 30, 2021}

\author[bsc]{Abel Gargallo-Peir\'o\corref{cor1}}
\ead{abel.gargallo@bsc.es}

\author[bsc]{Matias Avila}
\ead{matias.avila@bsc.es}

\author[bsc]{Arnau Folch}
\ead{arnau.folch@bsc.es}

\address[bsc]{Computer Applications in Science and Engineering, Barcelona Supercomputing Center, 08034 Barcelona, Spain.}

\cortext[cor1]{Corresponding author}

%
%
%

\begin{abstract}
A new topography adapted mesh generation process tailored to simulate Atmospheric Boundary Layer (ABL) flows on complex terrains is presented.
 The mesher is fully automatic given: the maximum and minimum surface mesh size, the size of the first element of the boundary layer, the maximum size in the boundary layer and the size at the top of the domain. The following contributions to the meshing workflow for ABL flow simulation are performed. First, we present a smooth topography modeling to query first and second order geometry derivatives. Second, we propose  a new adaptive meshing procedure to discretize the topography  based on two different metrics. Third, the ABL mesher is presented, featuring both prisms and tetrahedra. We extrude the triangles of the adapted surface mesh, generating prisms that reproduce the Surface Boundary Layer. Then, the rest of the domain is meshed with an unstructured tetrahedral mesh. In addition, for both the surface and volume meshers we detail a hybrid quality optimization approach, analyzing its impact on the solver for high-complexity terrains. We analyze the convergence of the triangle adaptive approach, obtaining quadratic convergence to the geometry and reducing to one half the error for the same amount of degrees of freedom than without adaptivity and optimization. We also study the mesh convergence of our RANS solver, obtaining quadratic mesh convergence to the solution, and using a 30\% of the degrees of freedom while reducing a 20\% of the error of standard semi-structured approaches. Finally, we present the generated meshes and the simulation results for a complete complex topographic scenario.
\end{abstract}


\end{frontmatter}


\section{Introduction}
\label{sec:intro}

During the last decades, the advances in Computational Fluid Dynamics (CFD) techniques together with the increase of the available computational power have widened the engineering applications that require of numerical analysis.
In wind resource assessment and wind power forecasting it is of special interest the simulation of Atmospheric Boundary Layer (ABL) flows. In the ABL, orographic gradients, ground surface drag, and atmospheric thermal instabilities from radiative forcing, can generate turbulence and strong wind shear (vertical velocity gradients) in the so-called Surface Boundary Layer (SBL), which extends up to a $10-20\%$ of the total ABL depth \cite{kopp1984remote,garratt1994atmospheric,blocken2007cfd,wizelius2007developing,stull2012introduction}. 
The requirement to reproduce the high-gradients of the SBL has been translated in most mesh generators for ABL simulation in a fully structured mesh in the normal direction to the surface \cite{sorensen:hypgrid,sorensen:ellypsys3Dvalidation,openfoam:web,gargallo:meshForABLandWindFarms,gargallo2018JCP:WindFarms,zephy}.
In addition to the structure along the vertical direction, most mesh generators specifically designed to discretize the ABL \cite{sorensen:hypgrid,sorensen:ellypsys3Dvalidation,openfoam:web,gargallo:meshForABLandWindFarms,gargallo2018JCP:WindFarms} are also semi-structured in the surface and do not feature adaptivity to the terrain. 

The use of structured grids was first used in Finite Difference and Finite Volumes applications, introducing the effect of the topography by means of a change of coordinates in the formulation of the problem instead of discretizing the topography with the mesh \cite{haltiner:numericalMeteorology,arakawa:vertical}.
Also in Finite Element and Finite Volume applications the use of structured hexahedral meshes represents the current most standard procedure. This  semi-structured meshing strategy exploits the advantage of hexahedra to align the mesh with the flow on offshore applications and also exploits the tensor structure of hexahedra to reproduce the close-to-surface boundary layer.
Several mesh generation approaches have based on solving a system of Partial Differential Equations locating the nodes to improve the orthogonality of the mesh and attain the desired element volume \cite{michelsen:1994block,sorensen:hypgrid,sorensen:ellypsys3Dvalidation,marras:atmosphericStabilization}.
In contrast with these approaches, in \cite{gargallo:meshForABLandWindFarms,gargallo2018JCP:WindFarms} it is proposed a procedure that combines sweeping the quadrilateral surface mesh with a quality-based mesh optimization, to determine the best configuration of the nodes of the mesh according to the chosen quality measure.

All those approaches feature hexahedra and use block structured strategies to provide a fine resolution in the zone of interest and less resolution far away from this region. However, since in each block the mesh is structured, inevitably finer resolution of the interest areas is extended to the rest of the domain. Although having the drawback of increasing the required number of elements of the mesh, this structured strategies have been exploited in offshore cases (analyzing the wind resource on the sea) or in topographic scenarios that do not feature high complexity  in order to align the mesh with the wind inflow direction, generating one mesh for each simulated wind inflow direction.
There are several alternatives to the use of hexahedral elements.
First, regarding the simulation of ABL flows, in \cite{zephy} a fully prismatic mesh is generated to discertize the ABL.
This work takes advantage of generating a triangle surface mesh to avoid extending the finer mesh size in the interest region to the rest of the domain and, simultaneously, takes advantage of the tensor direction of the prisms to discretize the boundary layer.
As an alternative to structured or semi-structured approaches, different approaches 
featuring tetrahedra
have been followed for problems that require 
a mesh
conformal with the topography but that do not require a boundary layer
close to the surface \cite{montero19983,montero2005genetic,escobar:smoothingSrfParam,escobar:tetMeshTerrains,behrens:amatos,oliver:windMeccano,gargallo:representingUrbanGeometries}.

Herein, we are interested in generating meshes in complex topographic scenarios, see for instance results in Section \ref{sec:exampleBadaia}. 
In those scenarios, even without Coriolis force, there is not a unique alignment direction close to the surface due to the interaction of the wind with the  topography. 
In addition, herein we model the Coriolis force in the flow, which translates in the fact that, even in  offshore scenarios, the wind direction close to the ground is different than the wind direction in the top of the domain. 
The change of the wind direction determined by Coriolis force is in average around 20 degrees from the ground to the top of the ABL, at 1.5km over the topography.
Taking into account both that the wind direction changes at different heights and also due to the influence of the topography, in this work we propose to use an unstructured mesh on the surface. 

Several contributions in mesh generation for Atmospheric Boundary Layer flow simulation are proposed in this work, each presented in a different section of the paper.
First, in Section \ref{sec:cfd} we present the target CFD model and we use it to detail the required meshing features.
Second,  we model the topography to allow queering first and second order derivatives to compute tangent and Hessian-based metrics, see Section \ref{sec:geometry}.
In addition, we incorporate the notion of metric complexity to define the Hessian-based metric accordingly to the desired number of degrees of freedom.
Third, we devise a new meshing approach for the surface topography based on a metric-driven mesh adaptation process, see Section \ref{sec:surface}.
Following, we discretize the ABL with a new hybrid approach, featuring prisms to resolve the SBL and tetrahedra to discretize the rest of the domain, see  Section \ref{sec:volume}.
Thus, out of the SBL the mesher is not constrained by the structure of the close to surface prismatic layer and unstructured tetrahedra are exploited to allow incresing the mesh size in all the directions.
Next, in Section \ref{sec:optimization} we present the quality and optimization framework for hybrid meshes developed in this work for topographic scenarios, analyzing the effect of the optimization for the target simulation.
Finally,   in Section \ref{sec:results}, we analyze the convergence to the geometry and to the solution of the proposed meshing approach and we compare it  to the standard strategy in the literature.
Together with obtaining the expected quadratic convergence, we also illustrate that with the proposed approach with a 30\% of the nodes of standard semi-structured strategies the error is reduced a 20\%.


\section{CFD model: Atmospheric Boundary Layer flow simulation}
\label{sec:cfd}
The meshing techniques developed in this work have been specifically designed to simulate Atmospheric Boundary Layer flows on complex topographies. 
Herein, we present a short overview of our CFD model for ABL flows, linking it with the required mesh specifications for our specific purpose mesher. 
Our complete CFD model, featuring the coupling with the actuator disc model for wind farm simulation, and the implementation, can be found in \cite{avila:CFDframworkWindFarms, Diaz2018}. 

\begin{figure}[!t] 
	\centering
	\newcommand{\xFigSize}{0.45\textwidth}
	\begin{tabular}{cc}
		\subfigure[]{\label{fig:escudo1}
			\includegraphics[width=\xFigSize]
			{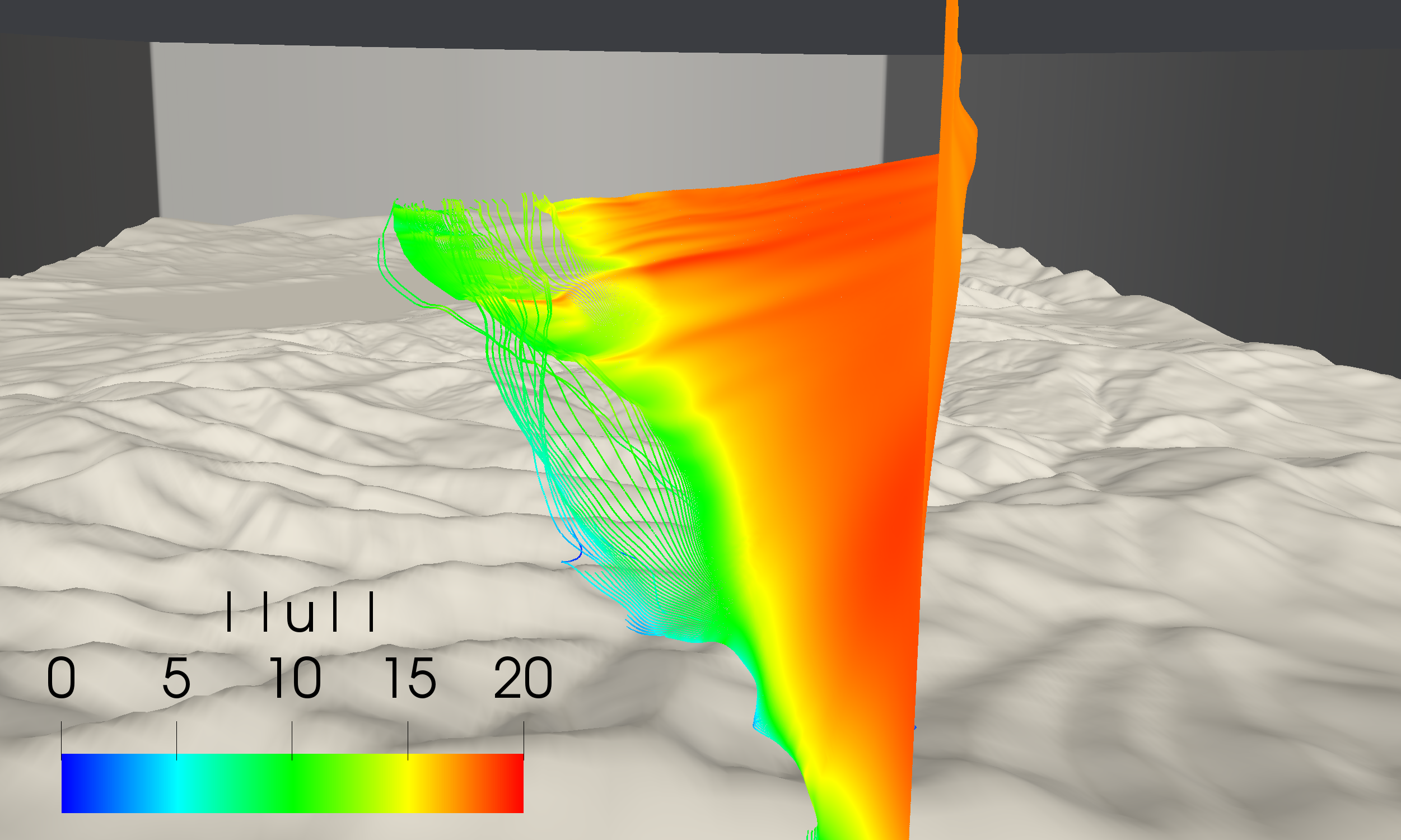}  	 }
		&
		\subfigure[]{\label{fig:escudo2}
			\includegraphics[width=\xFigSize]
			{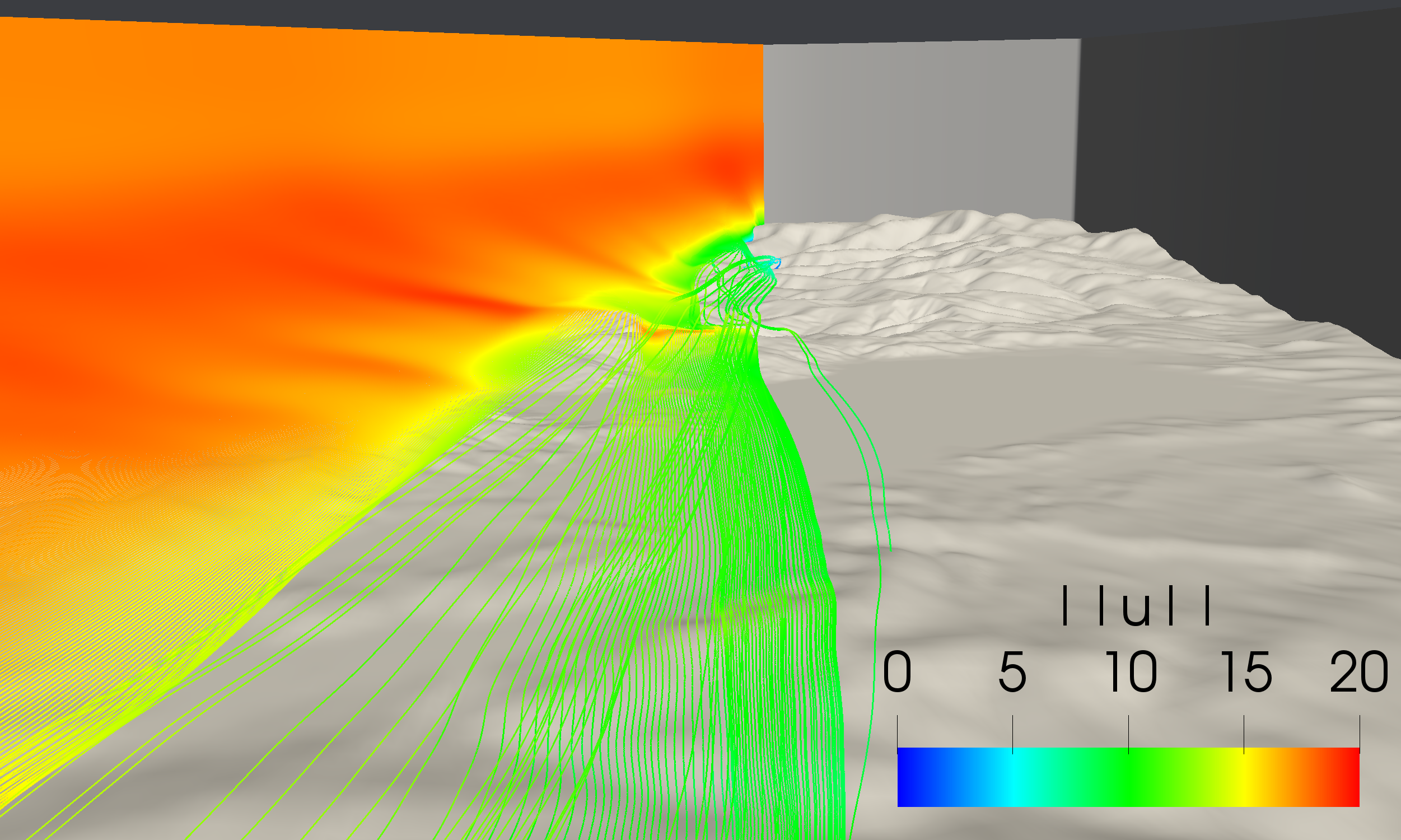}  	 } 	
	\end{tabular}
	\\
	\subfigure[]{\label{fig:escudo3}
		\includegraphics[width=0.7\textwidth]
		{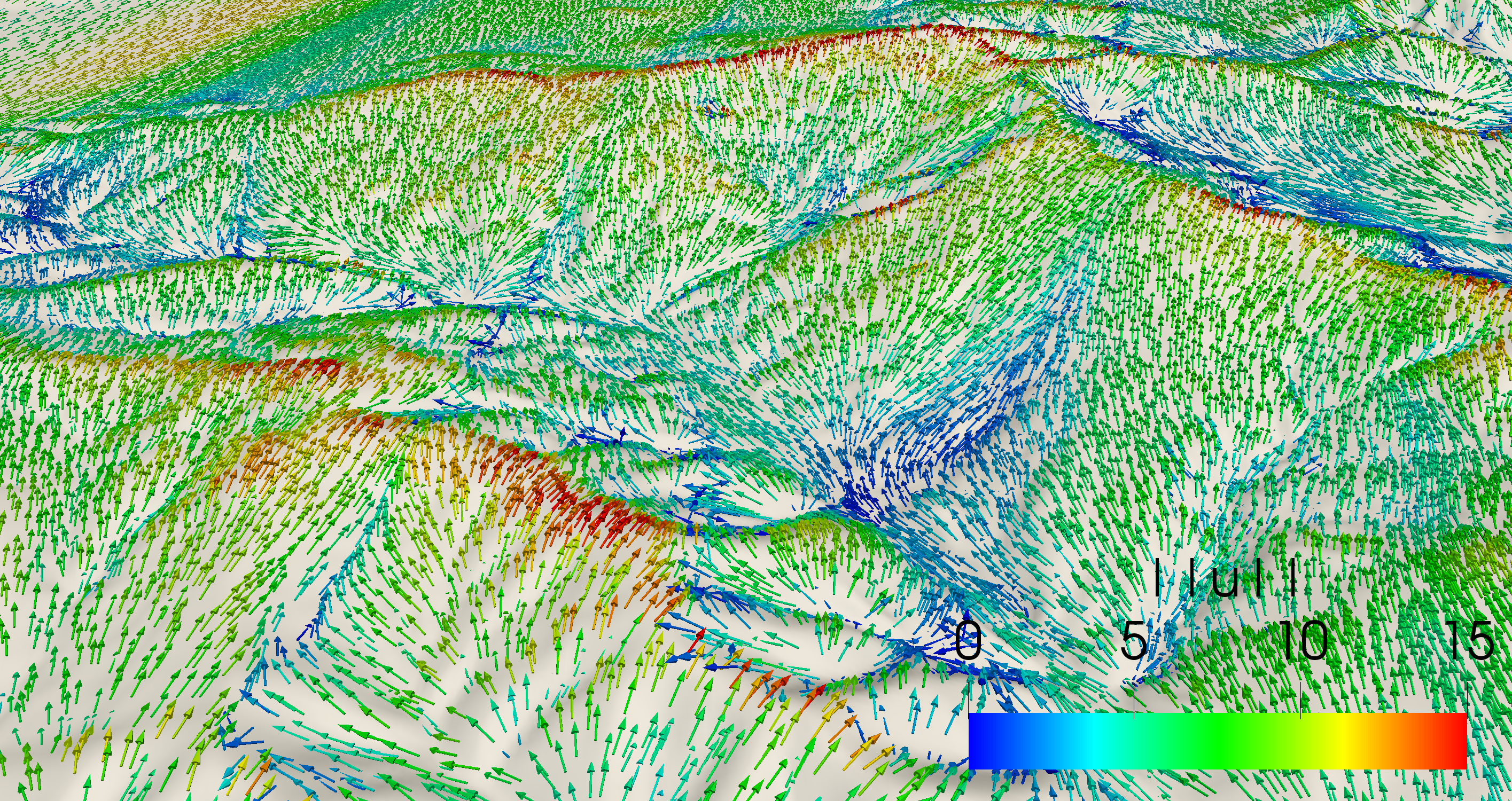}  	 } 	
	\caption{Simulation of the CFD model for the Escudo mountain range (Spain):
		\subref{fig:escudo1}-\subref{fig:escudo2}  perspective of the streamlines on an upstream vertical line viewed from the  inflow and outflow boundaries, respectively, and 
		\subref{fig:escudo3} velocity vector field at 5m from the ground. 
	}
	\label{fig:escudo}
\end{figure}

Considering the flow as incompressible and isothermal (neutral stability), the incompressible Navier Stokes equation coupled to the $k$-$\eps$ RANS model modified with the Apsley and Castro correction  for the mixing length limitation \cite{Apsley} in the presence of Coriolis forces is written as:
\begin{eqnarray}
	\nabla \cdot \mathbf{u} 
	&=& 0 \label{masseq} \\ 
	\frac{\partial \mathbf{u}}{\partial t} + \mathbf{u}\cdot \nabla \bu-\nabla \cdot 
	\left( \nu _t \nabla ^{s}\mathbf{u}\right)+  \nabla p + f_c\, \mathbf{ e}_z \times \left(\mathbf{u} -\mathbf u_g\right)  
	&=&
	\mathbf 0 \label{momeq} \\ 
	\frac{\partial k}{\partial t} +
	\bu\cdot \nabla k -\nabla \cdot \left(   
	\frac{\nu_t }{\sigma _{k}}  \nabla k\right)  
	+ \varepsilon - P_k  
	&=& 0  \label{keyeq} \\
	\frac{\partial \varepsilon }{\partial t} + \mathbf{u}\cdot \nabla
	\varepsilon -\nabla \cdot \left(\frac{\nu_t}{\sigma_\eps}
	\nabla \varepsilon \right) 
	-\frac{\varepsilon}{k}\left(C_{1}^\prime P_k - C_2\varepsilon  \right) &=&
	0  \label{epseq}
\end{eqnarray}
closed with the following expression for the turbulent viscosity
\begin{equation}
  \nu_t = C_\mu\frac{k^2}{\eps},
  \label{turvieq}
\end{equation}    
where the unknowns are the velocity field $\bu$, the  pressure $p$, the turbulent kinetic energy $k$, and the dissipation rate of turbulent kinetic energy
$\eps$. 
The term  $\mathbf u_g$ corresponds to the given geostrophic velocity, which is related to the geostrophic pressure gradient driving the flow, $\nabla p_g$, as $\nabla p_g=f_c\, \mathbf{ e}_z \times \mathbf u_g $.
We highlight that the fifth term on the left hand side (LHS) of momentum equation, Eq. \eqref{momeq}, is the Coriolis force, 
where the Coriolis parameter is $f_c=2\Omega \sin \lambda$  (with  $\Omega$ the earth's rotation rate and  $\lambda$ the latitude of the scenario) and $\mathbf{e}_z$ the local unit vector pointing in the vertical direction $z$. 
Regarding the mesh generation process, the Coriolis force term is responsible of not being possible to determine a unique horizontal alignment direction of the mesh with the flow, since it introduces a wind turning with height.
In Figures \ref{fig:escudo1} and \ref{fig:escudo2} we illustrate the simulation of the ABL flows with the CFD model on the Escudo mountain range in Spain.
In particular, we illustrate the streamlines (colored according to the wind speed) originated on a vertical line on the inflow boundary. It can be observed that the flow direction changes according to the point height.
Thus, the modeling of Coriolis force precludes the a priori alignment of the mesh with the flow without regards of the topographic scenario.

In the turbulence equations (\ref{keyeq})-(\ref{epseq}), the term $P_k=2\nu_t S$ is the kinetic energy production due to shear stress, with $S=\mathbf{\nabla}^s\bu:\mathbf{\nabla}^s\bu$ ($\mathbf \nabla ^s$ denotes the symmetric gradient operator). 
For the coefficients of the $k$-$\eps$ modified model ($C_\mu$, $C_2$, $\sigma_k$, $ \sigma_\eps$) we follow Panofsky and Dutton \cite{Panofsky},
and the  coefficient $C_1^\prime$  in the RHS of the dissipation equation (\ref{epseq}) corresponds to the modified coefficient accounting for Coriolis force, originally proposed by Apsley and Castro \cite{Apsley}. 

Proper boundary conditions need to be added to the Navier Stokes  (\ref{masseq})-(\ref{momeq}) and turbulence  $k$-$\eps$ (\ref{keyeq})-(\ref{epseq}) equations. The boundaries of the computational domain are split into inflow, outflow, bottom and top.
On the inflow boundary wind velocity $\bu$ and turbulence unknowns $k$ and $\eps$ are imposed as height profiles. These vertical profiles are generated from a single-column (1D) precursor simulation ({\it i.e.} that assumes flat terrain and uniform roughness).
Analyzing both the physics and the chosen precursor simulation, we can estimate the close to ground anisotropy in the SBL and use this a priori known data to set the mesher.
Herein, we use by default a first element height of 1 meter, a growth factor of the boundary layer is chosen in the interval $[1.05,1.2]$. 
On the outflow  boundary geostrophic pressure and zero shear stress are imposed for the momentum equation; while zero gradient in the normal direction is imposed to the turbulence unknowns.

On the top boundary, symmetry boundary conditions are imposed for velocity, that is, zero gradient for the tangential velocity component, and the normal velocity component  fixed to zero ({\it i.e.} $\bu\cdot\mathbf n =0$).  Zero normal gradient is  imposed to the turbulence unknowns. 
On the bottom boundary a wall law satisfying the Monin-Obukhov  \cite{monin1954basic} is imposed to the momentum and turbulence equations removing a boundary layer of thickness $\delta_w$.

As highlighted in Section \ref{sec:intro},
in complex topographies there is also not a unique flow direction due to the effect of the topography. Thus, in Fig. \ref{fig:escudo3} we plot 
the velocity field at 5m over the ground for the Escudo mountain range.
It can be observed that the wind is aligned with valleys and influenced by the different topography features and does not have a unique alignment.
Thus, the advantages of standard meshing approaches, featuring a structured strategy aligned with a main inflow direction, are precluded by both the complex topographic features and the simulation of the Coriolis force on the ABL flow.

We highlight that the model 
is implemented in the finite-element multi-physics parallel solver {\alya} \cite{houzeaux:alyaFractional,vazquez:alyaExascale}. However, the proposed meshing algorithm has been implemented in the external model-independent pre-process code \xWindMesh\ \cite{avila:CFDframworkWindFarms,gargallo:meshForABLandWindFarms,gargallo2018JCP:WindFarms,gargallo2017subdividing}. As a result, meshes generated with this utility can be used to simulate both with RANS or Large-Eddy Simulation (LES) turbulence models and, in addition, are also valid for solvers based on other numerical methods such as Finite Volumes.

\section{Topography geometry: modeling and metrics}
\label{sec:geometry}

In this section, we detail the modeling of the topography geometry performed in this work and the computation of the target metrics for mesh adaptation.
Herein, the topography geometry is defined with a piece-wise linear parameterization.
However, the computation of the target metrics will require to compute first an second order derivatives. Thus, first in Section \ref{sec:HOapproximant} we 
 model the topography to obtain a smooth parameterization where first and second order derivatives can be queried.
Second, in Section \ref{sec:topoMetrics} we present the two metrics that will be the driver to adapt the parametric mesh to capture the topograpy features. These metrics will be particularized for the case of topographies, and the notion of metric complexity will be used to resolve the topographic features with approximately the target number of degrees of freedom.

\subsection{Modeling of the topography for  first and second order derivative computation} 
\label{sec:HOapproximant}

This work is devoted to generate meshes conformal with the terrain.  
The geometry corresponds to real measured data that can be given in many formats, such as contour topography maps,  Cartesian grids or point clouds.
Herein, all the input frameworks are converted into a triangle mesh that is used as a geometry representation. 
Following, we define the proposed parameterization of the target surface (topography).
Due to the existent data extraction procedures, topographic data for wind resource assessment
 ensures that each point in the plane has a unique height value.
This is so since the triangle mesh that represents the geometry can be interpreted as a planar triangle mesh with a height value in each of the points of the mesh.
In particular, we define $\xParamPlane$ as the parametric plane,  the region in $\xreal{2}$ where the topography is defined. 
Thus, given a point $\xParamPoint{}\in\xParamPlane$ the height function of this point $\zh(\xParamPoint{})$ on the topography is defined as
\begin{equation}
\begin{array}{rccc}
\zh:&\xParamPlane\subset\xreal{2}&\longrightarrow&\xreal{}\\
&\xParamPoint{}=\xParamCoord{}&\longmapsto & \zh(\xParamPoint{}),
\end{array}
\label{eq:zh}
\end{equation}
where $\zh(\xParamPoint{})$ is computed by finding the triangle to which this point belongs and interpolating the height value in the triangle.
The function $\zh$ unequivocally determines a parameterization $\xparametrizationDiscrete$ of the topography surface $\xSurface$ as 
\begin{equation}
\begin{array}{rccc}
\xparametrizationDiscrete:&\xParamPlane\subset\xreal{2}&\longrightarrow&\Sigma\subset\xreal{3}\\
&\xParamPoint{}=\xParamCoord{}&\longmapsto & \xPhysicalCoord{}=(\xPhysicalCoordFirst{},\xPhysicalCoordSecond{},\zh(\xParamPoint{}))\xTransposed.
\end{array}
\label{eq:parameterizationDiscrete}
\end{equation}
To compute for instance the curvature of the target surface (Section \ref{sec:topoMetrics}) in the mesh adaptation process (Section \ref{sec:surface}), it is required to perform queries of the first and second order derivatives of the geometry. However, since the topography representation is a piece-wise linear triangle geometry, the derivatives of the topography geometry are not well defined. 
On the one hand, there are discontinuities of the first order derivative in the edges between the elements of the topography geometry. On the other hand, since the  geometry mesh features piece-wise linear elements, the curvature of each triangle configuring the geometry is null.
Thus, similarly to \cite{frey2000surface,frey2003surface},
herein we propose to reconstruct local high-order approximations of the geometry that allow queries of the derivatives of the surface representation.

To generate the high-order approximation in a point $\xParamPoint{}\in\xParamPlane$,  first it is located the triangle that contains this point in the  mesh that defines the geometry. Following, 
several layers of elements adjacent to the container elements are gathered.
Each new layer is computed as the triangles neighboring to the current considered elements. Then, we compute the least squares approximation of the desired order of the cloud of points determined by the neighborhood of elements around the target point. Specifically, given a set of topography points $\lbrace\xPhysicalCoord{1},\ldots,\xPhysicalCoord{n_p}\rbrace$ the high-order approximation used in this work is written as
\begin{equation}
\zp(x,y) = 
\sum_{\{i,j\}\in I_p} a_{ij}x^i y^j,
\label{eq:zp}
\end{equation}
where  
$\xOrderHO$ is the desired polynomial degree, 
$I_p=\{  \{i,j\} \ |\ i,j\ge 0 \text{ and } i+j\leq \xOrderHO\}$ is the set of indexes lower or equal to $\xOrderHO$,
and $\{a_{ij}\}_{\{i,j\}\in I_p}$ are the coefficients of the polynomial  $\zp(x,y) $ on $x$ and $y$.
In particular, we seek the approximation that better fits the cloud of points in the least squares sense:
\begin{equation}
\label{eq:leasSquaresZ}
\{a_{ij}\}_{\{i,j\}\in I_p}=
\argmin_{
	\substack{
		a_{ij}\in\xreal{} \\
		\{i,j\}\in I_p}}
\sum_{k=1}^{n_p} 
\left(
\zp(x_k,y_k) - z_k 
\right)^2
,
\end{equation}
where $\xPhysicalCoord{k}=(x_k,y_k,z_k)$ are the cloud point coordinates.

By default, $\xOrderHO$ levels of elements around the container element are considered. Once all the neighboring points are gathered, 
if the number of parameters is larger than the number of points,
an extra layer of neighboring elements is included, repeating this process until the number of points ensures a well-posed minimization problem. In particular, if we have the exact same number of parameters we would obtain a surface containing all the points of the cloud.

Once computed the smooth high-order approximation, the new parameterization $\xparametrizationP$ is defined as
\begin{equation}
\begin{array}{rccc}
\xparametrizationP:&\xParamPlane\subset\xreal{2}&\longrightarrow&\xreal{3}\\
&\xParamPoint{}=\xParamCoord{}&\longmapsto & \xPhysicalCoord{}=(\xPhysicalCoordFirst{},\xPhysicalCoordSecond{},\zp(x,y))\xTransposed.
\end{array}
\label{eq:parameterizationHO}
\end{equation}
Note that this parameterization has well defined first and second order derivatives that will allow to compute the desired metrics on the geometry in Section \ref{sec:topoMetrics}.
The order of approximation used in this work in the presented examples is three, and, accordingly to what has been previously stated, the number of levels of adjacency to compute the approximation is also three (without requiring in any example to automatically computing extra layers of neighbors).

\subsection{Metrics of the topography geometry}
\label{sec:topoMetrics}

Herein, we target a surface mesh that has the desired edge length and that reproduces the curvature of the topography.
Hence, the surface mesh will be generated by means of ensuring that the edges of the parametric mesh are unitary with respect to two different metrics: the metric of the tangent space, and the metric derived from the Hessian of the parameterization.
Following Section \ref{sec:HOapproximant}, we denote a general parameterization of a surface $\xSurface$ as:
\begin{eqnarray}
\label{eq:parameterizationGeneral}
\begin{aligned}
\xparametrization:\xParamPlane\subset\xreal{2}
&\longrightarrow\xSurface\subset\xreal{3}\\
\xParamPoint{}=\xParamCoord{}
&\longmapsto\xPhysicalCoord{}
= \xparametrization(\xParamPoint{}).
\end{aligned}
\end{eqnarray}
We define the curve  $\xcurve(t)$ between two points $\xPhysicalCoord{1}=\xparametrization(\xParamPoint{1})$ and $\xPhysicalCoord{2}=\xparametrization(\xParamPoint{2})$ on $\xSurface$ in terms of the parameterization $\xparametrization$ and the edge  $[\xParamPoint{1}, \xParamPoint{2} ]$ on the parametric space as
\begin{equation}
\label{eq:curveSurface}
\begin{aligned}
\xcurve:[0,1]&\longrightarrow\xSurface\subset\xreal{3}\\
t&\longmapsto
\xparametrization(\xParamPoint{}(t)), 
\end{aligned}
\end{equation}
where $\xParamPoint{}(t)$ is defined as
\begin{equation}
\label{eq:edgeParametricSpace}
\begin{aligned}
\xParamPoint{}:[0,1]&\longrightarrow[\xParamPoint{1}, \xParamPoint{2} ]\subset\xParamPlane\subset\xreal{2}\\
t&\longmapsto \xParamPoint{1}+t(\xParamPoint{2}-\xParamPoint{1}),
\end{aligned}
\end{equation}
and
with  $\xcurve$ fulfilling that $\xcurve(0) = \xPhysicalCoord{1}$ and $\xcurve(1)=\xPhysicalCoord{2}$. 
Following the ideas in \cite{Frey2008Book,gargallo20:topoAdapted} the length of the curve on the surface 
in terms of the parametric coordinates of the two surface nodes is:
\begin{eqnarray}
\label{eq:lengthCurveParamPoints}
\begin{aligned}
l_{\xMetricFisrtFF}(\xParamPoint{1},\xParamPoint{2}) 
&
= \int_0^1 
\left(
(\xParamPoint{2}-\xParamPoint{1}) \xTransposed
\ 
\xMetricFisrtFF(\xParamPoint{}(t)) 
\ 
(\xParamPoint{2}-\xParamPoint{1})
\right)^{1/2}
\    \text{d} t,
\end{aligned}
\end{eqnarray}
where 
\begin{eqnarray}
\xMetricFisrtFF(\xParamPoint{}(t)) 
\equiv
\xMetricFisrtFF
=
\xgrad\xparametrization\xTransposed
\cdot
\xgrad\xparametrization
=
\left( \parcial{\xparametrization}{\xParamCoordFirst{}}\ \  \parcial{\xparametrization}{\xParamCoordSecond{}} \right)  \xTransposed
\cdot
\left( \parcial{\xparametrization}{\xParamCoordFirst{}} \ \  \parcial{\xparametrization}{\xParamCoordSecond{}} \right) ,
\label{eq:metric1FF}
\end{eqnarray}
is the matrix expression of the first fundamental form of the surface $\xSurface$ at the point on the parametric space $\xParamPoint{}(t)$, 
see Eq. \eqref{eq:edgeParametricSpace}.
In particular, for isotropic mesh generation with a desired length $h$ of an edge known in each region of the domain, we define the \emph{tangent metric} in terms of $\xMetricFisrtFF$ as
\begin{equation}
\label{eq:metricTangent}
\xMetricTangent
:=
\frac{1}{h^2}\xMetricFisrtFF,
\end{equation}
and the corresponding length measure as
\begin{equation}
\label{eq:lengthMeasureTangent}
 	l_{\xMetricTangent}(\xParamPoint{1},\xParamPoint{2}) 
:= 
	\int_0^1 
\left(
(\xParamPoint{2}-\xParamPoint{1}) \xTransposed
\cdot 
\xMetricTangent(\xParamPoint{}(t)) 
\cdot
(\xParamPoint{2}-\xParamPoint{1})
\right)^{1/2}
\    \text{d} t.
\end{equation} 
Note that we would ideally like the edges of the mesh to have measure 1 with the metric $\xMetricTangent$.
In this manner, in the adaptive procedure that will be presented in Section \ref{sec:surface}, Eq. \eqref{eq:lengthMeasureTangent} will be used to compute the length of the curve on the surface in terms of the coordinates of the nodes in the parametric space. Next, these elements with measure greater than one (with a safety factor) will be refined until all edges fulfill the desired metric.

As previously introduced,
we target that the surface mesh also reproduces the curvature of the surface. That is,  the second order derivatives of the surface will be also taken into account in the adaptive procedure. 
To do so, we explicitly exploit
that the parameterization $\xparametrization$ in Eq. \eqref{eq:parameterizationGeneral} for topographic geometries can be rewritten as $\xparametrizationP(\xParamCoordFirst{},\xParamCoordSecond{})=(\xParamCoordFirst{},\xParamCoordSecond{},\zp(\xParamCoordFirst{},\xParamCoordSecond{}))$, as detailed in Eq. \eqref{eq:zp} and Eq. \eqref{eq:parameterizationHO} from Section \ref{sec:HOapproximant}.
Thus, the parameterization can be also understood as a field $\zp$ over a 2D mesh on $\xParamCoordFirst{}$ and $\xParamCoordSecond{}$.
Hence, herein we propose to
adapt the mesh to capture the curvature of the geometry using techniques for 2D mesh adaptation to reduce the interpolation error of $\zp(\xParamCoordFirst{},\xParamCoordSecond{})$. 

Following the ideas presented in \cite{peraire:87,frey2005anisotropic}, from the Hessian $\xHessianNice$ of the topography at a point, 
\begin{equation*}
	\xHessianNice=\left(
	\begin{array}{cc}
		\frac{\partial^2 \zp}{\partial x\partial x}&\frac{\partial^2 \zp}{\partial x\partial y} \\
		\frac{\partial^2 \zp}{\partial y\partial x} & \frac{\partial^2 \zp}{\partial y\partial y}
	\end{array}
	\right),
\end{equation*}
the following  metric is defined:
\begin{equation}
\label{eq:metricCurvature}
\xMetricCurvature^\beta=\ \xMat{V}\ (\beta\ \xMat{D})\ \xMat{V}^t,
\end{equation}
where $\xMat{V}=(\vect{e_1},\vect{e_2})$ is the matrix composed by the eigenvectors $\vect{e_1}$ and $\vect{e_2}$ of $\xHessianNice$, $\xMat{D}= \text{diag}(|\lambda_1|,|\lambda_2|)$ a diagonal matrix with its absolute value of the eigenvalues, and $\beta>0$ a curvature discretization real parameter.
Similarly to \cite{loseille2009optimal,loseille2011continuous1,loseille2011continuous2},
we replace the geometric notion of $\beta$, which allows us to refine further or less the mesh according to the curvature, to relate the amount of accuracy with the desired number of degrees of freedom.
Our objective with respect to the standard structured approaches is to locate the degrees of freedom where they are required. Thus, we  use this parameter $\beta$ to design a metric that is resolved with less or equal degrees of freedom than the ones that we would use in a structured mesh.
That is, we aim to generate a mesh with a similar number of nodes, but with increased accuracy.

In \cite{loseille2011continuous1,loseille2011continuous2} the complexity $\xComplexity$ of a metric $\xMetric$ is defined as 
\begin{equation}
\label{eq:complexity}
\xComplexity(\xMetric) :=\int_{\Omega}\sqrt{ \det \xMetric(\xParamPoint{}) } \ d\xParamPoint.
\end{equation}
The complexity of a metric 
can be interpreted as the counterpart at the continuous level of the number of nodes of a mesh. In particular, in \cite{loseille2011continuous1} it is proven that the following relation holds:
\begin{equation}
\label{eq:nodesComplexity}
\xNumNodes = \alpha\ \xComplexity(\xMetric)
\end{equation}
where $\xNumNodes$ is the number of nodes of the required unitary mesh, and $\alpha$ a constant depending on the domain and the mesh. The value of $\alpha$ is theoretically 2, but can change depending on the discrepancy between an ideal unit mesh and the generated discrete mesh. 
In \cite{loseille2011continuous2} different examples are presented featuring highly anisotropic metrics, where in practice values of $\alpha$ in the range of $[1.5,3]$ are obtained. Herein, 
 we use the theoretical value $\alpha=2$ so that we can couple a priori the curvature metric with an estimation of the maximum number of nodes that we allow to capture it.
Taking this into account, and in accordance to \cite{loseille2011continuous1}, from any metric $\xMetric$ we can define a new one so that it has the desired complexity:
\begin{equation}
\label{eq:metric_complexity}
\xMetric_{\xNumNodes}:=\left( \frac{\xNumNodes}{\xComplexity(\xMetric)} \right)^{\frac{2}{d}} \xMetric,
\end{equation}
where $\xNumNodes$ is the estimated desired number of nodes, and where $d=2$ in this work.

Herein, let $\xNumNodes$ be the desired number of nodes for our surface mesh, and let 
$C:=\xComplexity(\xMetricCurvature^1)$ be the complexity of $\xMetricCurvature^1$. Then, we define 
\begin{equation}
\label{eq:betaStar}
\beta^*:={\xNumNodes}/(\alpha C)
\end{equation}
 and define the \emph{curvature metric} as:
\begin{equation}
\label{eq:metricCurvComplex}
\xMetricCurvatureComplexity:=\xMetricCurvature^{\beta^*}.
\end{equation}
We highlight that the metric $\xMetricCurvatureComplexity$  has the desired complexity:
\begin{eqnarray*}
\xComplexity(\xMetricCurvatureComplexity)
&=&\int_{\xParamPlane}\sqrt{ \det \xMetricCurvatureComplexity(\xParamPoint{}) } \ d\xParamPoint{}
=\int_{\xParamPlane}\sqrt{ \det \xMetricCurvature^{\beta^*}(\xParamPoint{}) } \ d\xParamPoint{}
\\
&=&\int_{\xParamPlane}\sqrt{ \det \beta^*\xMetricCurvature^{1}(\xParamPoint{}) } \ d\xParamPoint{}
=\int_{\xParamPlane}\sqrt{ (\beta^*)^2\det \xMetricCurvature^{1}(\xParamPoint{}) } \ d\xParamPoint{}\\
&=&
\beta^*\int_{\xParamPlane}\sqrt{\det \xMetricCurvature^{1}(\xParamPoint{}) } \ d\xParamPoint{}
 \overset{C:=\xComplexity(\xMetricCurvature^1)}{=}
 \beta^*C
\overset{ Eq. \eqref{eq:betaStar}}{=}\frac{\xNumNodes}{\alpha},
\end{eqnarray*}
and thus, as desired and according to  Eq. \eqref{eq:nodesComplexity} , 
the metric $\xMetricCurvatureComplexity$ is approximately resolved with $\xNumNodes$ nodes.

Following, similarly to Eq. \eqref{eq:lengthMeasureTangent},  the length of the edge $[\xParamPoint{1},\xParamPoint{2}]$ with respect to the  metric $\xMetricCurvatureComplexity$ is defined as
\begin{eqnarray}
\label{eq:lengthMeasureCurvature}
	l_{\xMetricCurvatureComplexity}(\xParamPoint{1},\xParamPoint{2}) 
:= 	\int_0^1 
\left(
(\xParamPoint{2}-\xParamPoint{1}) \xTransposed
\cdot 
\xMetricCurvatureComplexity(\xParamPoint{}(t)) 
\cdot
(\xParamPoint{2}-\xParamPoint{1})
\right)^{1/2}
\    \text{d} t,
\end{eqnarray}
where $\xParamPoint{}(t)$ is defined in Eq. \eqref{eq:edgeParametricSpace}.
The value of the estimated number of nodes $\xNumNodes$ will be detailed in Section \ref{sec:surface}. 

In particular, if the user has prescribed an element size of $h$, the edges of a mesh should ideally fulfill $\xlengthMT=1$ and $\xlengthMC= 1$.
However, these two conditions may not be achievable simultaneously since, for instance, the curvature of the geometry may demand $\xlengthMT<<1$.
Following, in Section \ref{sec:surface} it is detailed the adaptive mesh generation proposed in this work, based on the two presented metrics. First, it will be generated a coarse topography mesh, which will be locally refined it until no edges of the mesh have length greater than 1, or if the  length of the edges is lower than a minimum value set by the user.
With these two conditions, the edge length of the surface mesh is controlled and it is ensured that the mesh  reproduces the curvature of the geometry.
\color{black}


\section{Topography adapted surface mesh generation} 
\label{sec:surface}

In this section, we detail our surface mesh adaptation procedure for topographic geometries. 
The topography is parameterized according to Section \ref{sec:HOapproximant}, and the edge lengths of the mesh are measured according to the tangent and curvature metrics, Eq. \eqref{eq:lengthMeasureTangent} and Eq.  \eqref{eq:lengthMeasureCurvature} respectively, both detailed in Section \ref{sec:topoMetrics}.

\begin{figure}[!t] 
	
	\begin{center}
		\includegraphics[width=0.45\textwidth]
		{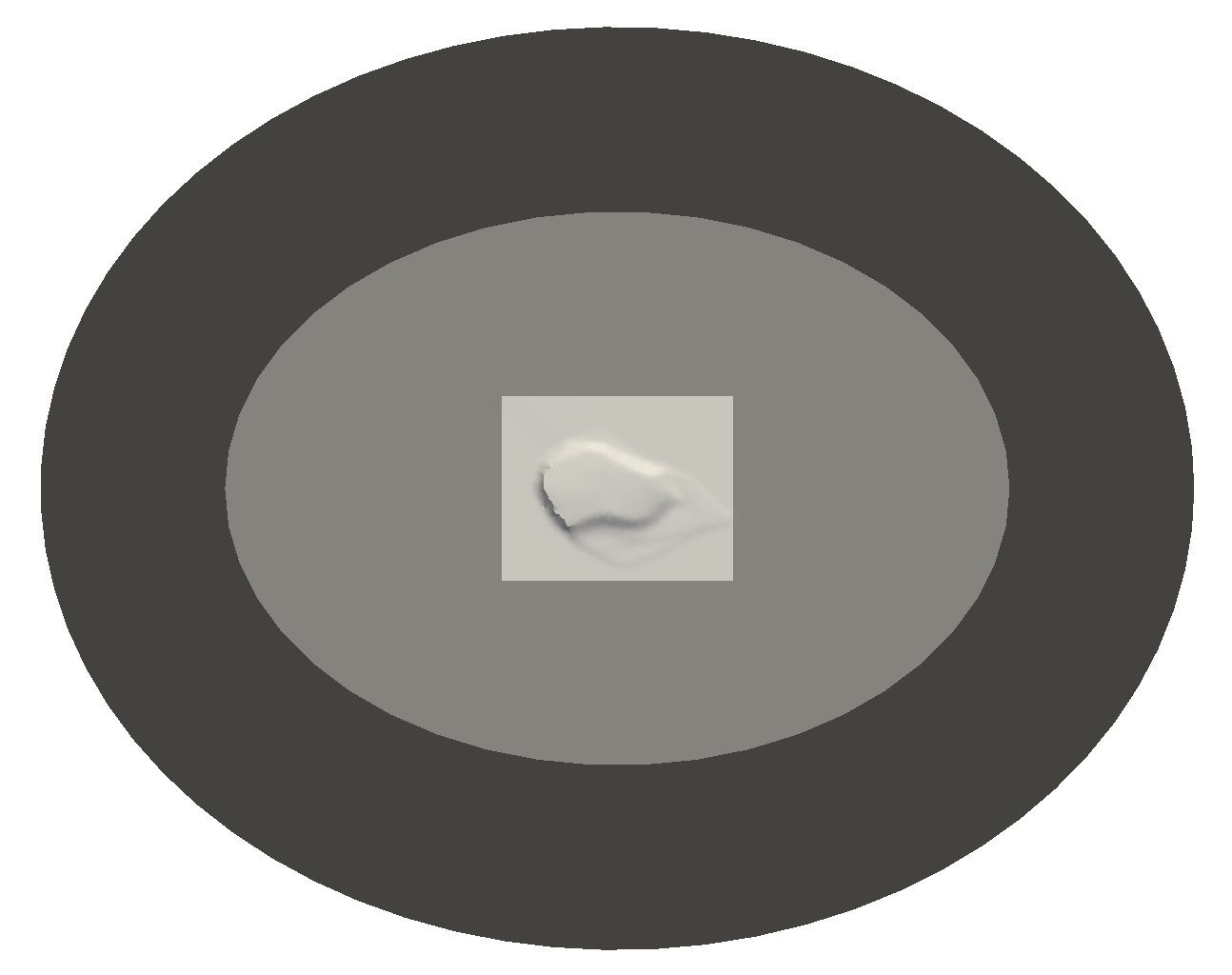}  	
	\end{center}
	
	\caption{Synthetic example illustrating the domain regions: farm (light gray), transition (gray) and buffer (dark grey).}
	\label{fig:windMeshRegions}
\end{figure}

In this work, the topography is divided into three regions, with three different levels of resolution illustrated in Figure \ref{fig:windMeshRegions}: the interest  farm area (higher resolution, light gray), a transition area (gray), and an elliptical buffer area to impose the boundary conditions (lower resolution, dark gray).
First, the farm area is meshed. The farm is a quadrilateral domain featuring the region of interest in the simulation (for instance, area where a wind farm is to be designed). 
The adaptive process to the topography is applied in the farm region, where  higher geometric accuracy is required to discretize the features of the topography.
The transition area  is an elliptical domain that encircles the farm region and is meshed with a triangle mesh that  smoothly matches the fine mesh of the transition with the element size of the buffer, the outer region. 
Finally, it is defined an additional elliptic region with a coarse element size  to impose the boundary conditions.

\newcommand{\xFigSize}{0.3\textwidth}
\begin{figure*} 
	\centering
	\begin{tabular}{ccc}
		\subfigure[]{\label{fig:triBolund_0}
			\includegraphics[width=\xFigSize]
			{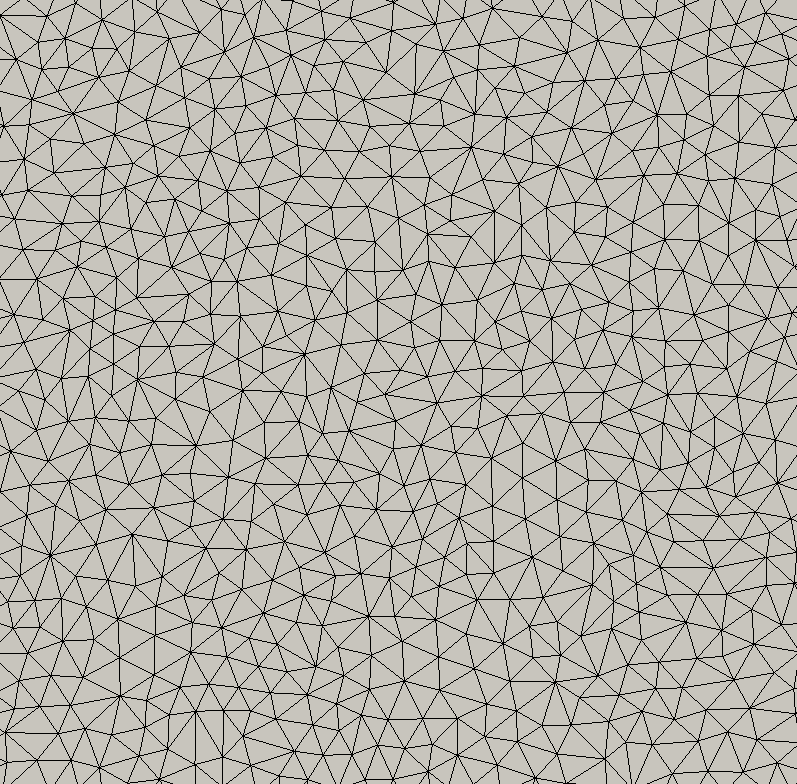}  	 }
		&
		\subfigure[]{\label{fig:triBolund_1}
			\includegraphics[width=\xFigSize]
			{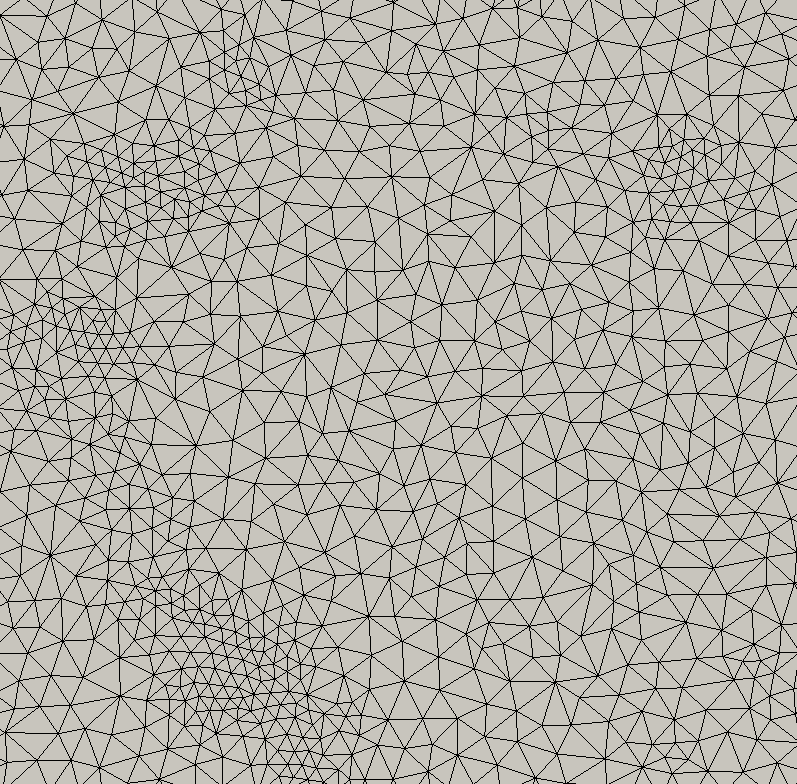}  	 } 	
		&
		\subfigure[]{\label{fig:triBolund_2}
			\includegraphics[width=\xFigSize]
			{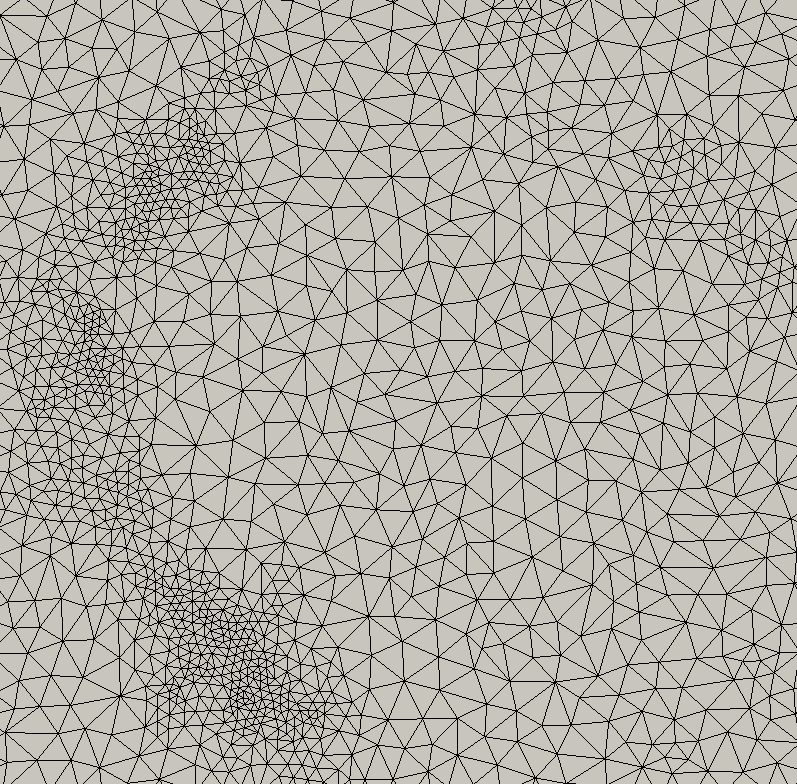}  	 }
		\\
		\subfigure[]{\label{fig:triBolund_3}
			\includegraphics[width=\xFigSize]
			{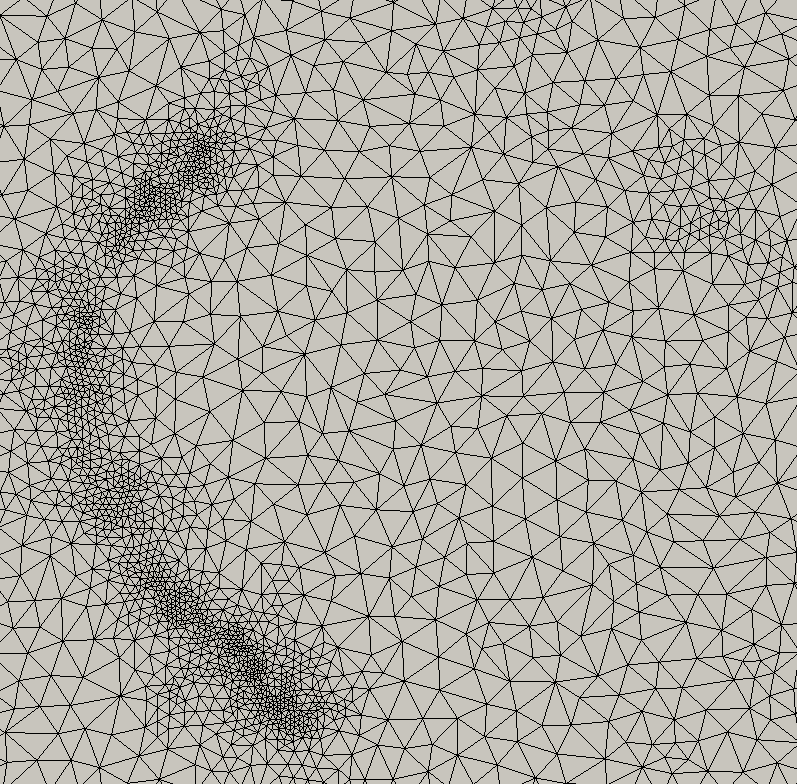}  	 } 	
		&
		\subfigure[]{\label{fig:triBolund_4}
			\includegraphics[width=\xFigSize]
			{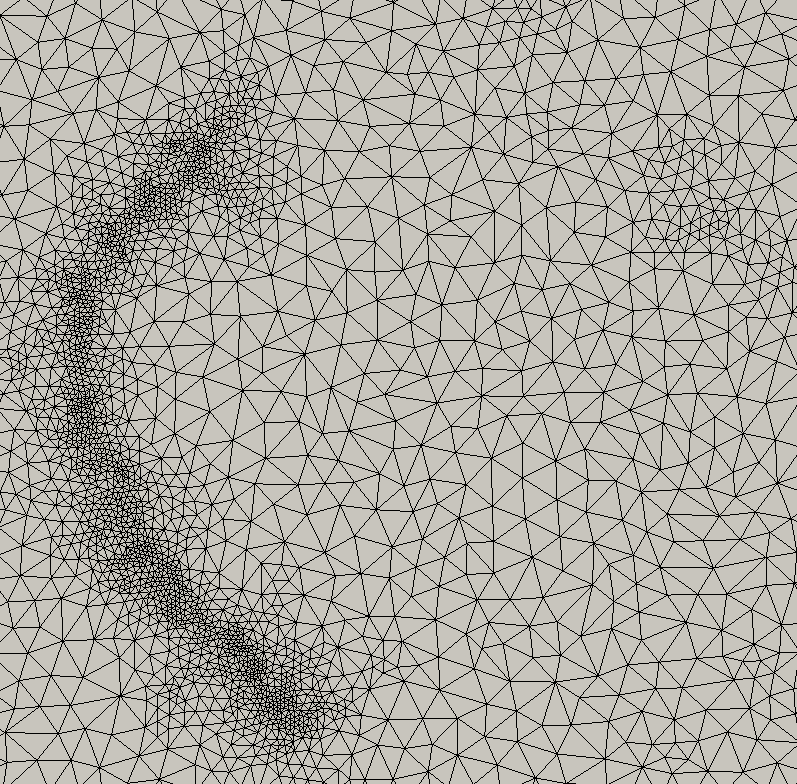}  	 }
		&
		\subfigure[]{\label{fig:triBolund_5}
			\includegraphics[width=\xFigSize]
			{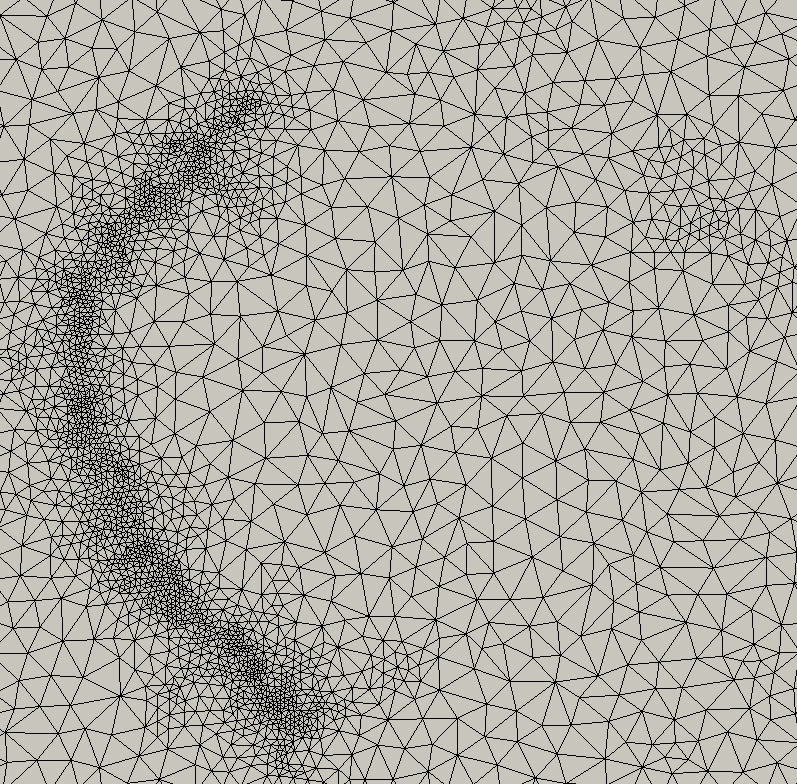}  	 } 	
		\\
		\subfigure[]{\label{fig:triBolund_6}
			\includegraphics[width=\xFigSize]
			{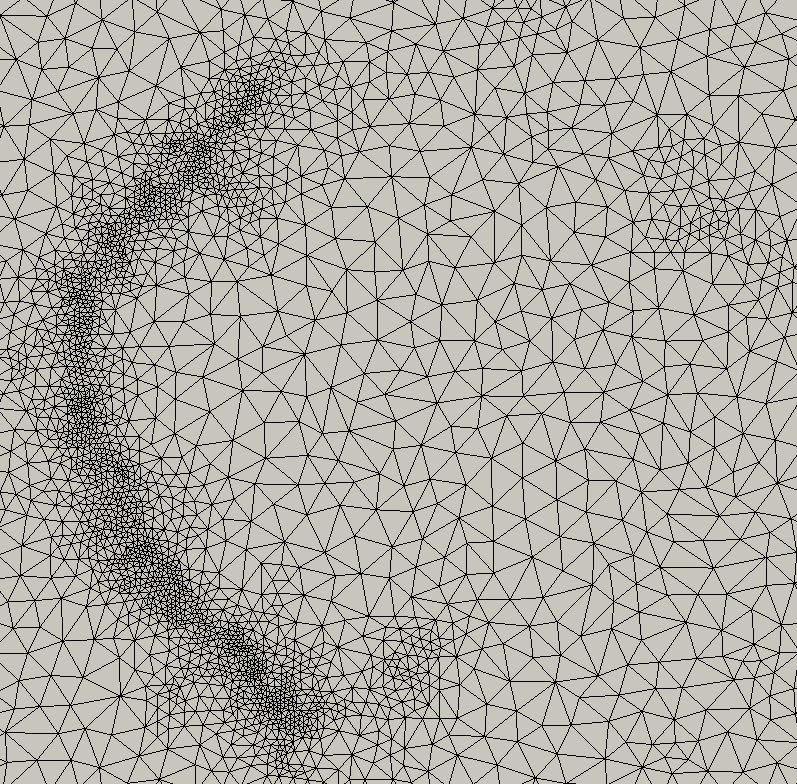}  	 }
		&
		\subfigure[]{\label{fig:triBolund_7}
			\includegraphics[width=\xFigSize]
			{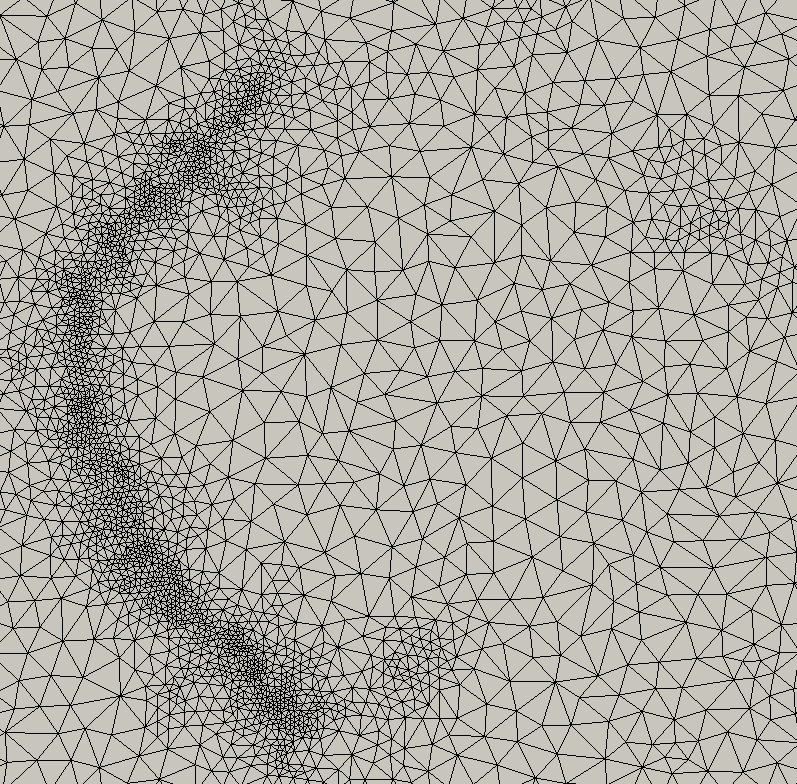}  	 }
		&
		\subfigure[]{\label{fig:triBolund_SRF}
			\includegraphics[width=\xFigSize]
			{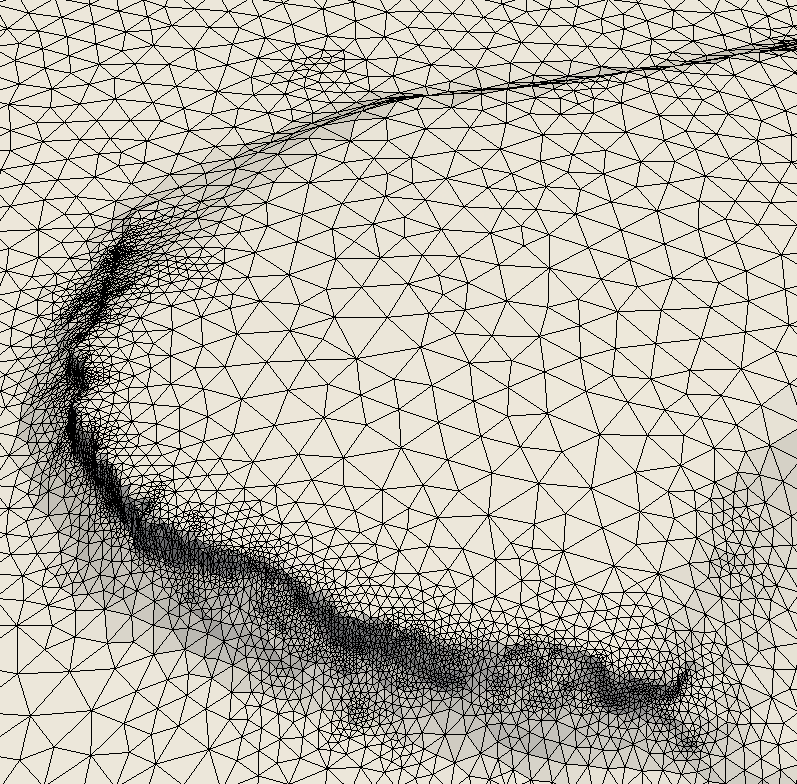}  	 } 
	\end{tabular}
	\caption{
		\subref{fig:triBolund_0}-\subref{fig:triBolund_7} Mesh adaptive process for the Bolund geometry illustrated in Figure \ref{fig:windMeshRegions}. (i) Final adapted topography surface mesh.
	}
	\label{fig:bolundExampleAdaptive}
\end{figure*}

To generate a surface mesh adapted to the topographic features in the interest region, two inputs are required: the maximum element size $h_{max}$ and the minimum allowed size $h_{min}$.
With these two inputs, we readjust the curvature metric complexity to allow resolving it without exceeding the used degrees of freedom.
In particular we define the maximum number of nodes to resolve the curvature metric, $\xNumNodes$ in Eq. \eqref{eq:betaStar} and Eq. \eqref{eq:metricCurvComplex} as the number of nodes that would require to a uniform mesh to discretize the farm region with an average size form the maximum $h_{max}$ and the minimum $h_{min}$ allowed sizes.

The adaptive mesh generation process is an iterative procedure  detailed in Algorithm \ref{alg:adaptiveSurface}. 
It starts by generating
a planar triangle mesh with the maximum element size, Line \ref{state:initialMesh},  generated using the Triangle mesh generator \cite{shewchuk1996triangle}. 
Following, in Line \ref{state:elemsToRefine}, the edge lengths of all the elements of the mesh are measured using the tangent and curvature metrics. If an edge of an element is greater than 1, this implies that the edge is longer than desired with the metric. 
The condition to accept the element as it is or refine it is relaxed to those elements with lengths greater than $\sqrt{2}$, see \cite{frey2000surface,Frey2008Book}. 
Thus, this element is included in a list of elements to be refined. 

Once all the mesh elements have been checked, we locally refine the mesh where it is required, see Line \ref{state:refineMesh2}. 
The refinement is performed by using the previous mesh as background mesh and asking to those elements included in the refinement list to have half of the size of the triangle.
This process is repeated until all the elements have length below the desired threshold for both metrics, or if the minimum edge length is below the desired minimum length. We highlight that the minimum length is checked using the Euclidean metric to control the minimum length in the mesh for the simulation.

\begin{algorithm}[t!]
	\caption{Surface mesh generation process adapted to the topography}
	\label{alg:adaptiveSurface}
	\begin{algorithmic}[1]
		\renewcommand{\algorithmicrequire}{\textbf{Input:}}
		\renewcommand{\algorithmicensure}{\textbf{Output:}}
		\Require{Topography surface $\xSurface$, Maximum edge length $h_{max}$, Minimum edge length $h_{min}$}
		\Ensure{Mesh $\xMesh$}
		\Function{MeshTopography}{$\xSurface$,$h_{max}$,$h_{min}$}
		\State $\xMesh$ $\gets$ generate planar triangle mesh of size $h_{max}$
		\label{state:initialMesh}
		\State \listElemsToRefine $\gets$  findElemsToRefine($\xMesh$,$\xSurface$,$h_{min}$) 
		\label{state:elemsToRefine}
		\While {\listElemsToRefine $\neq\emptyset$}
		\label{state:while}
		\State $\xMesh$ $\gets$ refinePlanarMesh($\xMesh$,\listElemsToRefine)
		\label{state:refineMesh2}
		\State \listElemsToRefine$\gets$findElemsToRefine($\xMesh$,$\xSurface$,$h_{min}$) 
		\label{state:listElemsToRefine2}
		\EndWhile
		\State  $\xMesh$ $\gets$ mapToSurface($\xMesh$,$\xSurface$)
		\label{state:mapToTopo}
		\State  $\xMesh$ $\gets$ optimizeMeshOnSurface($\xMesh$,$\xSurface$)
		\label{state:optimizeOnTopo}
		\State \Return  $\xMesh$
		\EndFunction
	\end{algorithmic}
\end{algorithm}

Once the adaptive process is finalized, it is obtained  a mesh that has elements of the desired size on the surface, and that reproduces the curvature of the geometry up to the minimum allowed mesh size and up to the metric complexity. 
Up to this point local mesh refinement has been performed in order to improve the accuracy of the geometric approximation.
However, in this process the quality of the generated mesh for simulation has still not been assessed. Thus, to conclude the generation of the topography surface mesh, in Line \ref{state:optimizeOnTopo} a quality optimization of the surface elements is performed.
The quality measures used in this work, together with the complete surface optimization procedure are detailed in Section \ref{sec:optimization}.

In Figure \ref{fig:bolundExampleAdaptive} the adaptive process is applied to generate a mesh on the Bolund peninsula (Denmark) geometry.
In this illustrative example, the input mesh sizes are $h_{max}=5$m and $h_{min}=0.5$m.
First, in Figure \ref{fig:triBolund_0} a planar mesh of constant element size is generated, Line \ref{state:initialMesh} in Algorithm \ref{alg:adaptiveSurface}. Then Figures  \ref{fig:triBolund_1} to \ref{fig:triBolund_7} illustrate the successive refinements of the mesh according to Line \ref{state:refineMesh2}. Finally, in Figure \ref{fig:triBolund_SRF} the final surface mesh is illustrated. In the procedure, the initial number of elements is 5089, and the final one is 9037. 
It can be observed that nodes have been located in the slopes and curved features of the terrain. Away from the features of the geometry, the mesh size is the  $h_{max}$, whereas close to the geometry features fine refinement has been performed complying the already detailed constrains.


\section{Hybrid Atmospheric Boundary layer mesh generation}
\label{sec:volume}

The Atmospheric Boundary Layer mesh is generated in a domain enclosed by the topography, a planar ceiling at the desired height (user input defaulted as 2km over the highest topography point), and an elliptic lateral wall. The elliptic lateral boundary is defined by extruding vertically the boundary of the 2D domain, see Figure \ref{fig:windMeshRegions}, up to the desired height. 

This volumetric domain is meshed following Algorithm \ref{alg:hybridABL}. First, in Line \ref{state:prism}, the triangle surface mesh of the topography is swept (extruded) to generate a structured prismatic mesh close to the ground.  
Once the prismatic mesh is finished, in Line \ref{state:tet}
an unstructured tetrahedral mesh is generated to fill the rest of the domain with tetrahedra: from the last prismatic layer to the flat ceiling located at height $z_{top}$.
Finally, the hybrid mesh is optimized to improve the quality of the generated mesh, Line \ref{state:optimizeHybrid}.

\begin{algorithm}[t!]
	\caption{Generation of the volumetric hybrid ABL mesh}
	\label{alg:hybridABL}
	\begin{algorithmic}[1]
		\renewcommand{\algorithmicrequire}{\textbf{Input:}}
		\renewcommand{\algorithmicensure}{\textbf{Output:}}
		\Require{Surface mesh $\xMesh_\xSurface$, Initial Size SBL $h_{0}$, Final Size SBL $h_{1}$, Growing Ratio $r$, Height SBL $z_{BL}$, Top Ceiling $z_{top}$}, Top Ceiling Size $h_{2}$
		\Ensure{Mesh $\xMesh$}
		\Function{MeshABL}{$\xMesh_\xSurface$,$h_{0}$,$h_{1}$,$h_{2}$,$r$,$z_{BL}$,$z_{top}$}
		\State $\xMesh_{\mathrm{pri}}\gets$ {MeshSBL}({$\xMesh_\xSurface$},$h_{0}$,$h_{1}$,$r$,$z_{BL}$)
		\label{state:prism}
		\State $\xMesh_{\mathrm{tet}}\gets$ GenerateTetMesh($\xMesh_{\mathrm{pri}}$,$z_{top}$,$h_{2}$)
		\label{state:tet}
		\State $\xMesh\gets$ GenerateHybridMesh($\xMesh_{\mathrm{pri}}$,$\xMesh_{\mathrm{tet}}$)
		\label{state:hybrid}
		\State $\xMesh\gets$ OptimizeQuality($\xMesh$)
		\label{state:optimizeHybrid}
		\EndFunction
	\end{algorithmic}
\end{algorithm}

Algorithm \ref{alg:prismaticMesh} details the prismatic generation process.
As input for the prismatic meshing process, it is required to provide the initial height of the elements ($h_{0}$), the growing ratio ($r$),  a maximum elemental height for the boundary layer ($h_{1}$), and the height until which the boundary layer extends ($z_{BL}$).
Herein, according to the CFD model described in Section \ref{sec:cfd}, the growth factor of the boundary layer is chosen in the interval $[1.05,1.2]$, and the anisotropy in the first layer is of the order of $1/100$, depending on the CFD case and the region of the domain.
In addition, the maximum height of the structured boundary layer region is set by default at the $20\%$ of the total ABL (400 meters over 2km, SBL), see \cite{kopp1984remote,garratt1994atmospheric,blocken2007cfd,wizelius2007developing,stull2012introduction}. If the elements become isotropic before reaching the height of the structured region, the sweeping is continued 
keeping constant the size along the vertical direction.

Given a initial triangle surface mesh generated with the approach presented in Section \ref{sec:surface}, the volume mesh of the close to the surface region is generated by means of an iterative sweeping procedure that requires two main steps to compute each new sweeping layer.
First, Line \ref{alg:sweepTriangleMesh} of Alg. \ref{alg:prismaticMesh}, given a layer of triangle elements, a new layer of triangles is generated.
These two triangle meshes determine the new prismatic layer by connecting each node on the initial layer with the corresponding swept node in the extruded layer, Line \ref{alg:generatePrisms}.
Finally, Line \ref{alg:optimizePrisms}, once a new layer of prisms has been generated, 
it is performed a local optimization of the low quality elements that have been generated up to this point.

\begin{algorithm}[t!]
	\caption{Generation of the prismatic mesh of the SBL }
	\label{alg:prismaticMesh}
	\begin{algorithmic}[1]
		\renewcommand{\algorithmicrequire}{\textbf{Input:}}
		\renewcommand{\algorithmicensure}{\textbf{Output:}}
		\Require{Surface mesh $\xMesh_\xSurface$, Initial height $h_{0}$, Final height $h_{1}$, Growing Ratio $r$, Height SBL $z_{BL}$}
		\Ensure{Mesh $\xMesh$}
		\Function{MeshSBL}{{$\xMesh_\xSurface$},$h_{0}$,$h_{1}$,$r$,$z_{BL}$}
		\State $\xMesh^{tri}_0\gets\xMesh_\xSurface$
		\State $\xMesh^{pri}_{}\ \gets\emptyset$
		\State $n\ \gets0$
		\While{min height of $\xMesh^{tri}_n$ $< z_{BL}$} 
		\State $n\ \gets n+1$
		\State $\xMesh^{tri}_n\ \gets$sweepTriangleMesh($\xMesh^{tri}_{n-1}$,$n$,$r$,$h_0$,$h_1$)
		\label{alg:sweepTriangleMesh}
		\State $\xMesh^{pri}_{n}\ \gets$ generatePrisms($\xMesh^{tri}_{n-1},\xMesh^{tri}_n$)
		\label{alg:generatePrisms}
		\State $\xMesh^{pri}\ \gets\xMesh^{pri}\cup\xMesh^{pri}_{n}$
		\label{alg:generatePrismMesh}
		\State $\xMesh^{pri}_{}\ \gets$ optimizePrismaticMesh($\xMesh^{pri}_{}$)
		\label{alg:optimizePrisms}
		\EndWhile
		\EndFunction
	\end{algorithmic}
\end{algorithm}


Given a layer of triangle elements, we generate a new layer of prisms by means of sweeping each node computing a new extruding length and an extruding direction. 
The current extrusion length is computed in a standard manner using a geometrical law based on a user input growing ratio, herein defaulted as 1.15.
As input for the user it is also asked a maximum element size (defaulted as the isotropic configuration of the swept element). Once this size is reached or if the element becomes isotropic, the extruding length is kept constant and not further increased with the growing ratio.

To determine the extruding direction, we compute the pseudo-normal \cite{roca:PhDDissertation,roca2010automatic} of the nodes adjacent to each swept node.  In particular, given a node $\xnode{}$ with $\xNumNodesLoop$ neighboring nodes $\{\xnode{1},\ldots,\xnode{\xNumNodesLoop}\}$, the pseudo-normal $\xPseudoNormal$ is defined as
\begin{equation*} 
	\begin{aligned}
		\xPseudoNormal:=&\frac{\xNewSum_{i=1}^{\xNumNodesLoop} \ \xnode{i}\times\xnode{i+1}}{\xNewNormSpace{\xNewSum_{i=1}^{\xNumNodesLoop} \ \xnode{i}\times\xnode{i+1}}{}}
		= 
		\frac{\xNewSum_{i=1}^{\xNumNodesLoop} \ (\xnode{i}-\xnode{})\times(\xnode{i+1}-\xnode{})}{\xNewNormSpace{\xNewSum_{i=1}^{\xNumNodesLoop} \ (\xnode{i}-\xnode{})\times(\xnode{i+1}-\xnode{})}{}}
		,
	\end{aligned}
\end{equation*} 
where we consider $\xnode{\xNumNodesLoop+1}\equiv\xnode{1}$.
The pseudo-normal is used to maximize the orthogonality of the extruded layer with respect to the previous one.
The main property of the pseudo-normal $\xPseudoNormal$ is that it defines the plane that maximizes the area of the projection of the polygon defined by $\{\xnode{1},\ldots,\xnode{\xNumNodesLoop}\}$, see \cite{roca:PhDDissertation,roca2010automatic}, reducing the chances of generating invalid extruded elements.
We highlight that the pseudo-normal is blended with the vertical direction in order to enforce that the mesh grows towards the ceiling and that it gets to the top orthogonally to the planar ceiling. 

\begin{figure}[!t] 
	\centering
	\hspace{-0.55cm}
	
	\begin{tabular}{cc}
		\subfigure[]{\label{fig:prismMesh}
			\includegraphics[width=0.4\textwidth]
			{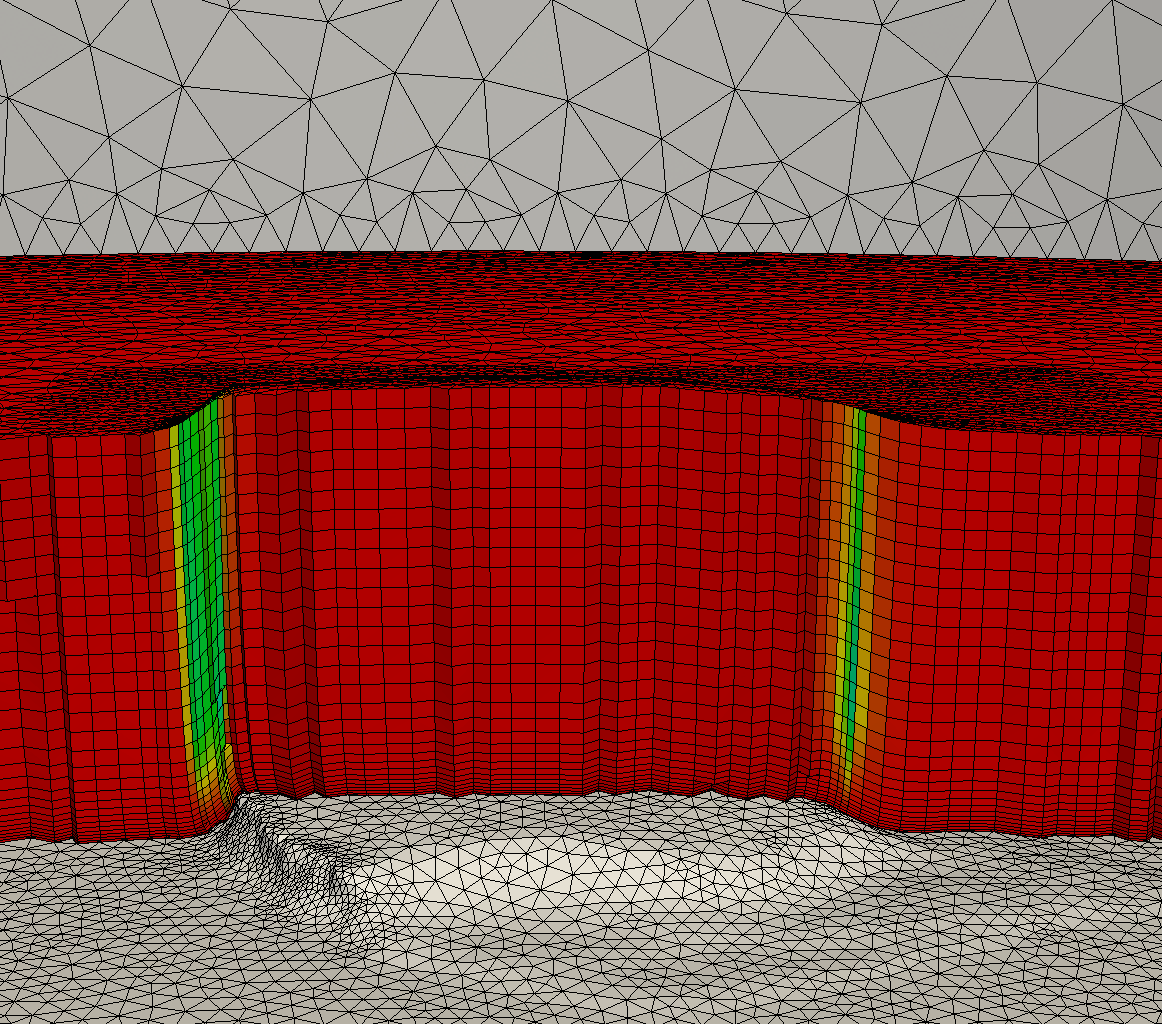}  	 } 	
		& 
		\subfigure[]{\label{fig:hybridMesh}
			\includegraphics[width=0.4\textwidth]
			{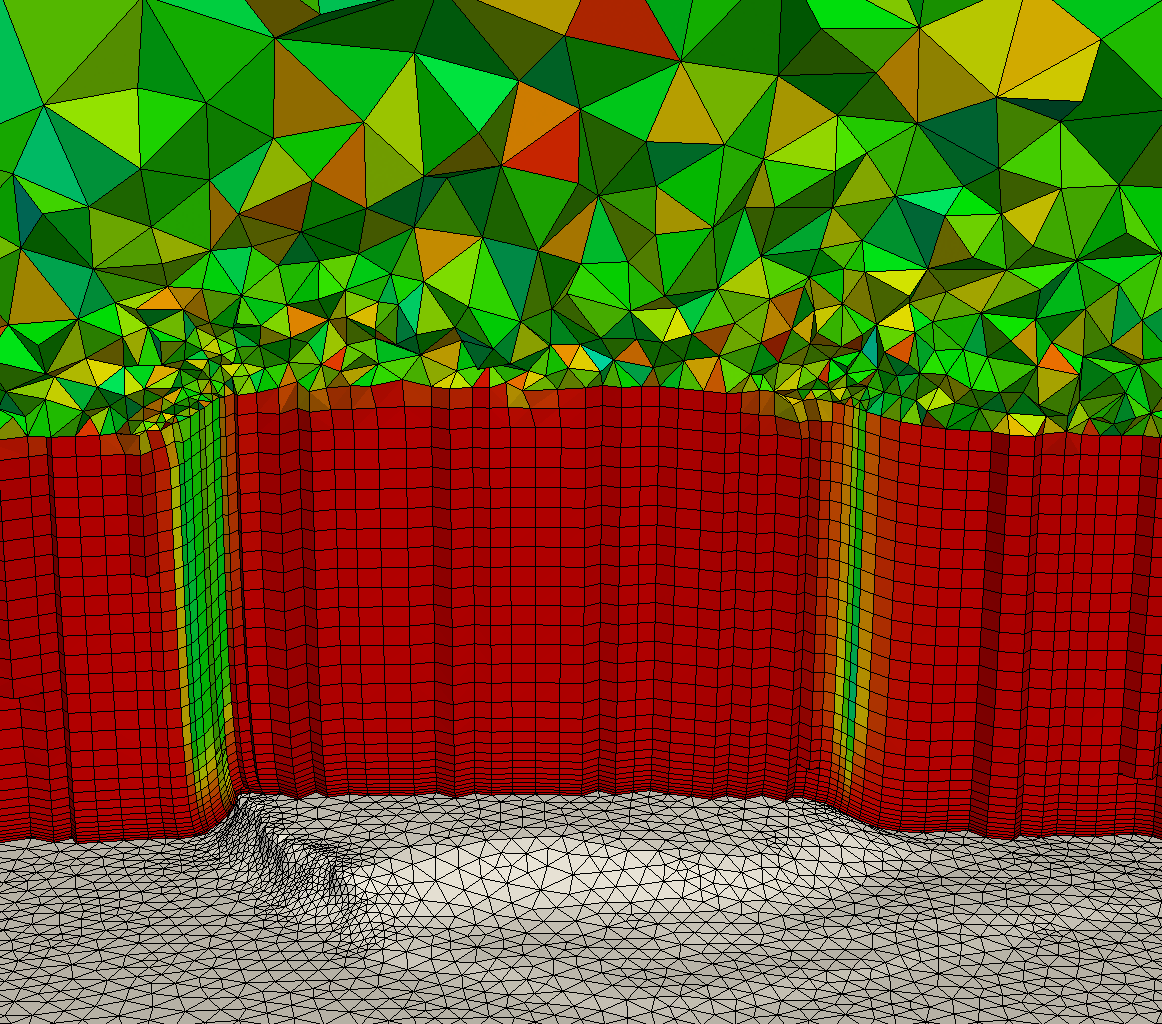}  	 } 
	\end{tabular}
	\includegraphics[width=0.2\textwidth]{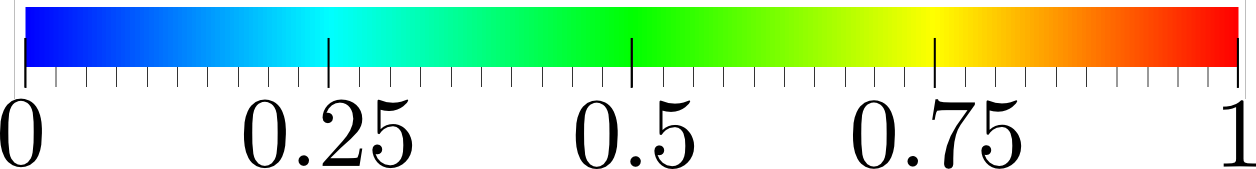} 
	\caption{
		Main meshing steps to generate the mesh of the Atmospheric Boundary Layer for the Bolund hill: \subref{fig:prismMesh} prismatic mesh and \subref{fig:hybridMesh} optimized hybrid prismatic-tetrahedral mesh.
		The elements are colored according to their quality. 
	}
	\label{fig:bolundWindMesh}
\end{figure}

It must be remarked that although the use of the pseudo-normal reduces the possibility of obtaining inverted or low quality elements, 
the complexity of the terrain may still drive to low quality configurations.
Thus, it is required to assess the quality of the mesh during the procedure to ensure its validity and, if possible, improve its quality.
The optimization process is detailed in Section \ref{sec:hybridOptimization}.
Once a new prismatic layer is generated, it is optimized to improve its configuration. 
Therefore, before generating a new layer of elements,  the current layer is optimized. 
Figure \ref{fig:bolundWindMesh}  illustrates the generated hybrid volume mesh. 
Figure \ref{fig:prismMesh} shows the first step of the process, the prismatic mesh obtaining by sweeping the topography surface mesh. 

Once the sweeping process is finalized, the unstructured tetrahedral mesher \tetgen\ \cite{si:tetgen} is used to generate a volume mesh that fills the rest of the domain up to the desired height, Line \ref{state:hybrid} of Alg. \ref{alg:hybridABL}. To do so, the target volume is defined in terms of a planar ceiling, an elliptical surface that encloses all the perimeter of the initial surface mesh, and the top boundary of the prismatic mesh. Each triangle of the boundary that encloses the volume is assigned with an element size field that determines the size of the tetrahedral mesh.
In particular, the triangles from the swept prismatic mesh are assigned with their own size.
The triangles of the planar top ceiling are assigned the element size chosen by the user, $h_2$ in Alg. \ref{alg:hybridABL}.
Finally, the triangles from the elliptic lateral are assigned an element size computed using a linear blend between the size at the swept surface and the size at the planar top ceiling.

After the topology of the hybrid mesh has been set, we perform a final  mesh optimization, see Section \ref{sec:hybridOptimization}, to compute the coordinates of the mesh nodes that deliver minimal distortion, according to the measure detailed in Section \ref{sec:qualityHybrid}.
In Figure \ref{fig:hybridMesh} we present the resulting  hybrid mesh for the Bolund hill.
We highlight that the domain specifications (size of the domain, mesh size and height of the boundary layer) have been set
to clearly illustrate in a unique figure the adaptation to the terrain, the prismatic boundary layer and the tetrahedral mesh.


\section{Hybrid mesh quality and optimization}
\label{sec:optimization}

In this section, we present the framework of quality measures for any element used in this work to assess the validity of the mesh. 
Following the ideas for high-order elements presented in \cite{gargallo:generation3Doptimization,AGPXRNJPGJSR:IMR13ewc}, we use the families of Jacobian-based quality measures to write the quality of an element in terms of the Jacobian of the mapping between an ideal element and the physical one. Thus, the quality of any element type can be assessed with the same functional if we are able to determine an ideal configuration.
First, in Section \ref{sec:qualityHybrid} we detail the framework that we use in this work for triangles, prisms and tetrahedra.
Next, in Section \ref{sec:optimization} we detail the volumetric quality optimization for hybrid prismatic-tetrahedral meshes to improve the quality of the generated ABL meshes.
In Section \ref{sec:optimizationSurface} we present the quality optimization of the surface meshes constrained to the topography.
Finally, in Section \ref{sec:optimizationResults} we analyze the improvements on the mesh derived from both the surface and volume optimization.

\subsection{Quality measures for hybrid meshes}
\label{sec:qualityHybrid}

To measure if an element is valid, and to quantify how much it differs from the desired ideal configuration, we use a distortion measure (see  \cite{knupp:algebraicQualityInitialMeshes} for a review of measures). 
A distortion measure quantifies in the range $[1,\infty)$ the deviation of an element with respect to an ideal configuration  (for instance, the equilateral triangle with the desired size for the triangle case). 
In this work, the distortion of an element with nodes $\xnode{1},\ldots,\xnode{\numberNodesHO{p}}$ is denoted as $\xDistortionHO(\xnode{1},\ldots,\xnode{\numberNodesHO{p}})$. 
The distortion takes value 1 when the element presents the desired configuration, and tends to infinity when the element degenerates.  
Following the ideas for high-order elements in \cite{gargallo:generation3Doptimization,AGPXRNJPGJSR:IMR13ewc},
the distortion measure used in this work can be written for any given element with nodes $\xnode{1},\ldots,\xnode{\numberNodesHO{p}}$ as:
\begin{equation}
\xDistortionHO(\xnode{1},\ldots,\xnode{\numberNodesHO{p}})
:=
{
	\frac{
		\xNewNormSpace{
			\eta_{sh}(
			\xDifferential{\xrepresentationIP}
			)
		}{\xIdealElem{}}
	}
	{
		\xNewNormSpace{1}{\xIdealElem{}}
	}
},
\label{eq:distortionElementGeneral}
\end{equation}
where  $\eta_{sh}$ is defined in \cite{knupp:algebraicQualityInitialMeshes} as 
\begin{equation*}
	\eta_{sh}(
	\xDifferential{\xrepresentationIP}
	):= 
	\frac
	{\xNewNormSpace{\xDifferential{\xrepresentationIP} }{F}^2}
	{\xDim\ \mid{\det\xDifferential{\xrepresentationIP} }\mid^{2/\xDim}},
\end{equation*}
where
$\xrepresentationIP(\xnode{1},\ldots,\xnode{\numberNodesHO{p}})$ is the mapping between the ideal $ \xIdealElem{} $ and physical elements,
$\xDifferential{\xrepresentationIP}$ is its Jacobian,
$\xNewNormSpace{\cdot}{F}$ is the Frobenius norm,
$\xNewNormSpace{\cdot}{\xIdealElem{}}$ is the $\mathcal{L}^2$ norm on the ideal element, 
$\xNewNormSpace{1}{\xIdealElem{}}$ is the measure of the ideal element,
and $\xDim$ is the dimension of the element ($\xDim=2$ for planar and surface meshes, and $\xDim=3$ for volumetric meshes).
In particular, the quality of an element is the inverse of the distortion:
\begin{equation}
\label{eq:quality}
\xQuality{}:=\frac{1}{\xdeviationLinear{}}\in[0,1],
\end{equation}
being 0 an invalid configuration, and $1$ the desired one.


\subsection{Hybrid tetrahedral-prism optimization}
\label{sec:hybridOptimization}

In this section, we present the quality optimization framework for hybrid tetrahedral-prismatic meshes.
To raise an optimization framework, first it is required to compute the quality/distortion of the mesh elements.
As detailed in Section \ref{sec:qualityHybrid}, in order to determine the quality/distortion of an element, an ideal element is required to be set.
We  set the ideal of each prism 
as an orthogonal prism with the desired anisotropy.
In particular, the ideal  of a given physical prism is defined as the initial surface element extruded orthogonally with the computed extrusion length at the layer where the prism has been generated. 
This is indeed the element that we would desire to generate, since it is extruded orthogonally and has the best possible configuration 
on the surface, since it has been optimized for the topography following the approach that will be detailed in Sec. \ref{sec:optimizationSurface}.
For the tetrahedral case, the ideal element is the equilateral tetrahedron. 

Once we are capable of determining the quality of the mesh, an optimization of the mesh can be derived.
In particular, we will minimize the distortion of the mesh in the least-squares sense \cite{gargallo:generation3Doptimization}.
Herein, the mesh distortion is a functional of the nodes of the mesh defined as:
\begin{equation}
	\label{eq:leastSquares3D}
	f(\xnode{1},\ldots,\xnode{\xNumNodes})=
	\xExtraConstant 
	\xNewSum_{\xelind=1}^{\xNumElements}
	\xDistortionHOind{\xelind}^2(\xnode{{\xelind}_1},\ldots,\xnode{{\xelind}_\numberNodesHO{p}}),
\end{equation}
where $\xNumNodes$ is the number of nodes of the mesh, $\xNumElements$ is the number of elements in the mesh, $\xDistortionHOind{\xelind}$ is the distortion of the $\xelind$-th element, and  $\xelind_i$ corresponds to the global node id of the $ith$ node of element $\xelind$. 

The optimization of the mesh corresponds to finding the node coordinates $\{\xnode{i}^*\}_{i=1,\ldots,\xNumNodes}$ such that
\begin{equation}
\label{eq:meshMinimization3D}
\{\xnode{1},\ldots,\xnode{\xNumNodes}\}
=
\argmin_{ \{\xnode{i}\in\xreal{3}\}_{i=1,\ldots,\xNumNodes}}
f(\xnode{1},\ldots,\xnode{\xNumNodes}).
\end{equation}
The minimization in Eq. \eqref{eq:meshMinimization3D} is performed following a Gauss-Seidel approach, as detailed in \cite{gargallo:generation3Doptimization}. That is, given a node $\xnode{i}$, we fix the rest of nodes of the mesh and we solve the non-linear minimization:
\begin{equation}
	\label{eq:leastSquaresNode}
	\xnode{i} =
	\argmin_{\xnode{i}\in\xreal{3}}
	\xExtraConstant 
	\xNewSum_{\xelind=1}^{\xNumElementsNeigh{i}}
	\xDistortionHOind{\hat{\xelind}}^2(\xnode{{\xelind}_1},\ldots,\xnode{{\xelind}_\numberNodesHO{p}})
	,
\end{equation}
where $\{{\hat{\xelind}}\}_{e=1,\ldots,\xNumElementsNeigh{i}}$ denotes the set of adjacent elements to node $\xnode{i}$, $\hat{\xelind}$ denotes the global id of the ${\xelind}$-th adjacent element of node $\xnode{i}$, $\xDistortionHOind{\hat{\xelind}}$ is the distortion of the ${\xelind}$-th element, and  ${\xelind}_i$ corresponds to the global node id of the $i$th node of the ${\xelind}$th neighbor element.

In the mesh generation process, the mesh quality optimization is applied in two different frameworks. 
On the one hand, once the meshing procedure is finished and the topology of the mesh is determined, we 
optimize the mesh to improve the quality of the final generated elements, Eq. \eqref{eq:leastSquares3D}, see  Line \ref{state:optimizeHybrid} of Alg. \ref{alg:hybridABL}
On the other hand, in the sweeping process, Line \ref{state:prism} in Alg. \ref{alg:hybridABL}, we perform for each swept layer, a local optimization to improve the current generated configuration.
During the prismatic sweeping process, Line \ref{alg:optimizePrisms} in Alg. \ref{alg:prismaticMesh},  
we do not optimize all the nodes that have been generated up to the current layer.
Instead, we derive a process that is it computationally efficient but that allows generating valid meshes. 
In particular, we only optimize the elements 
with a quality that is below an acceptance threshold.
The quality threshold for this local optimization is defaulted in our applications as 0.2.
To increase the freedom of the nodes that determine those low quality elements to find higher quality configurations, we include several layers of adjacency of the low quality elements in the procedure.
That is, we  find the elements with quality lower than a threshold, and we redefine the functional in Eq. \eqref{eq:leastSquares3D} in terms of the low quality elements and several layers of neighbors.
The boundary nodes of the neighborhood, and the rest of nodes of the domain are kept fixed.
Herein, by default two levels of neighboring elements are included in the local optimization process.

\subsection{Constrained surface mesh optimization}
\label{sec:optimizationSurface}

In this section, we present the quality optimization framework of the surface mesh constrained to the topography.
Similarly to the case for volume elements, to set the optimization framework, first it is required to compute the quality/distortion of the surface elements.
For each element on the surface, its ideal is defined as an equilateral triangle of the desired size. 
Once, each physical element has an ideal, the surface mesh can be optimized to minimize the elemental distortion.
It has to be taken into account that this ideal configuration can not in general be achieved since the mesh topology is now fixed and since the surface elements have their nodes tied to the geometry.

In contrast to the volumetric optimization presented in Eq. \eqref{eq:meshMinimization3D},
the surface mesh optimization is a constrained problem written as:
\begin{equation}
\label{eq:meshMinimizationSurfaceConstrained}
\lbrace\xnode{1}^*,\ldots,\xnode{\xNumNodesSurface}^*\rbrace
=
\argmin_{\xnode{1},\ldots,\xnode{\xNumNodesSurface} \in \xSurface}
f(\xnode{1},\ldots,\xnode{\xNumNodesSurface}),
\end{equation}
where $\xSurface$ is the topography surface, $f$ is the mesh distortion presented in Eq. \eqref{eq:leastSquares3D}  applied to the triangle surface mesh, and $\xNumNodesSurface$ are the number of surface nodes.
Using the surface parameterization, any surface node $\xnode{}$ can be rewritten as $\xparametrization(\xParamPoint{})$. Thus, 
similarly to  \cite{gargallo2013ADMOS,gargallo:LinearSurfaces,gargallo:highorderSurfaces},
we rewrite Eq. \eqref{eq:meshMinimizationSurfaceConstrained} as an unconstrained problem expressed on the parametric space $\xParamPlane$ of the surface:
\begin{equation}
\label{eq:meshMinimizationSurface}
\lbrace\xParamPoint{1}^*,\ldots,\xParamPoint{\xNumNodesSurface}^*\rbrace
=
\argmin_{\xParamPoint{1},\ldots,\xParamPoint{\xNumNodesSurface} \in \xParamPlane}
f(\xparametrization(\xParamPoint{1}),\ldots,\xparametrization(\xParamPoint{\xNumNodesSurface})).
\end{equation}
That is, the initially constrained surface optimization problem has been translated to finding
the location of the nodes on the parametric plane $\xParamPlane$ such that provide minimal elemental distortion (maximum quality) of the surface elements in the least squares sense.
Similar to the volumetric case, Sec. \ref{sec:hybridOptimization}, the minimization problem in Eq. \eqref{eq:meshMinimizationSurface} is solved in a Gauss-Seidel approach, see Eq. \eqref{eq:leastSquaresNode}.
Once the optimization is finalized, the nodes are mapped to the surface with the parameterization $\xparametrization$.
For the optimization process, we use the input piece-wise linear parameterization in order to ensure that the final location of the mesh nodes is on the input topography.

\subsection{Analysis of the impact of optimization in the mesh generation process}
\label{sec:optimizationResults}

\begin{figure*}[!t] 
	\centering
	\hspace{-0.55cm}
	\begin{tabular}{cc}
		\subfigure[]{\label{fig:surfaceMesh1}
			\includegraphics[width=0.45\textwidth]
			{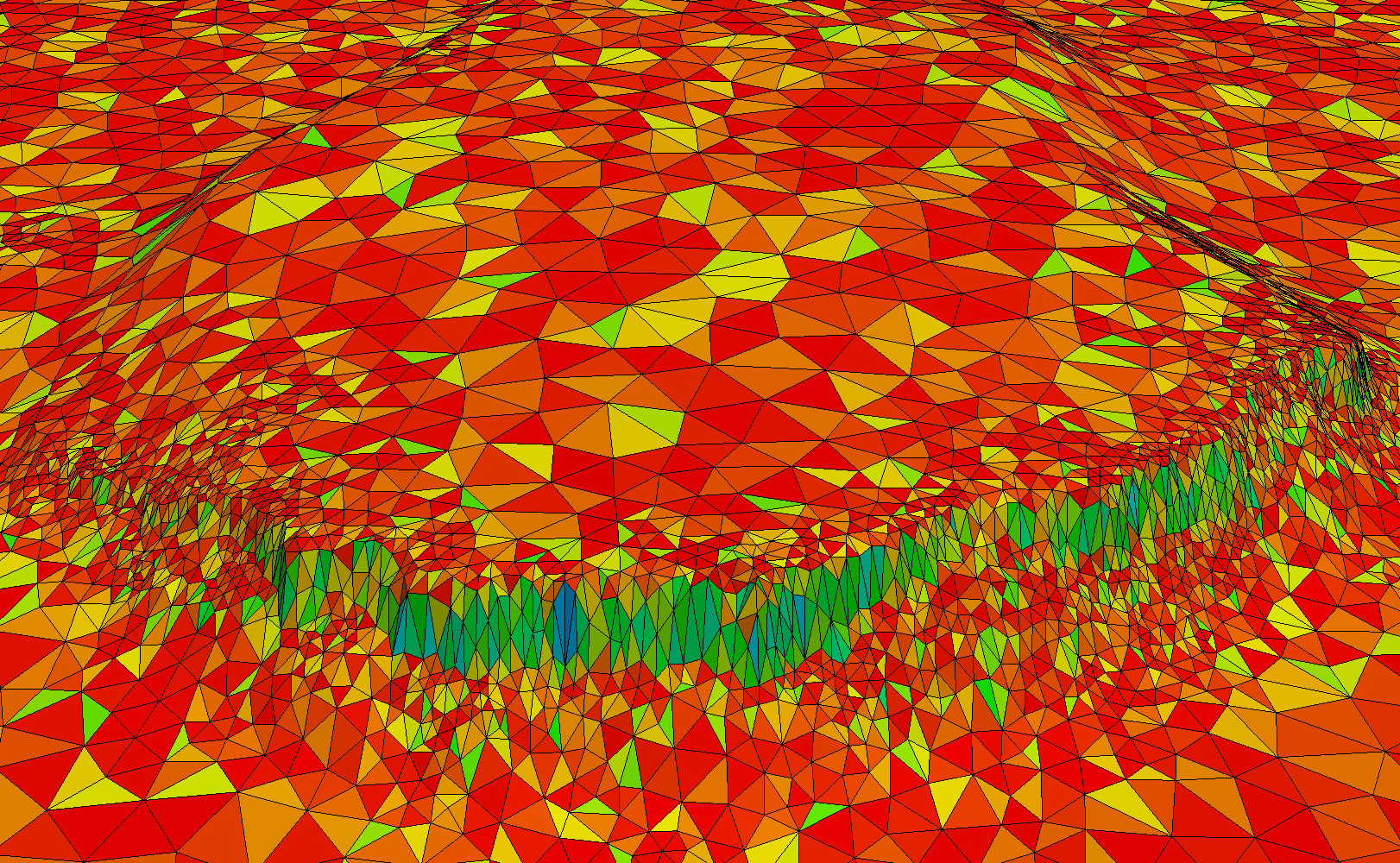}  	 } 	
		&
		\subfigure[]{\label{fig:surfaceMesh2}
			\includegraphics[width=0.45\textwidth]
			{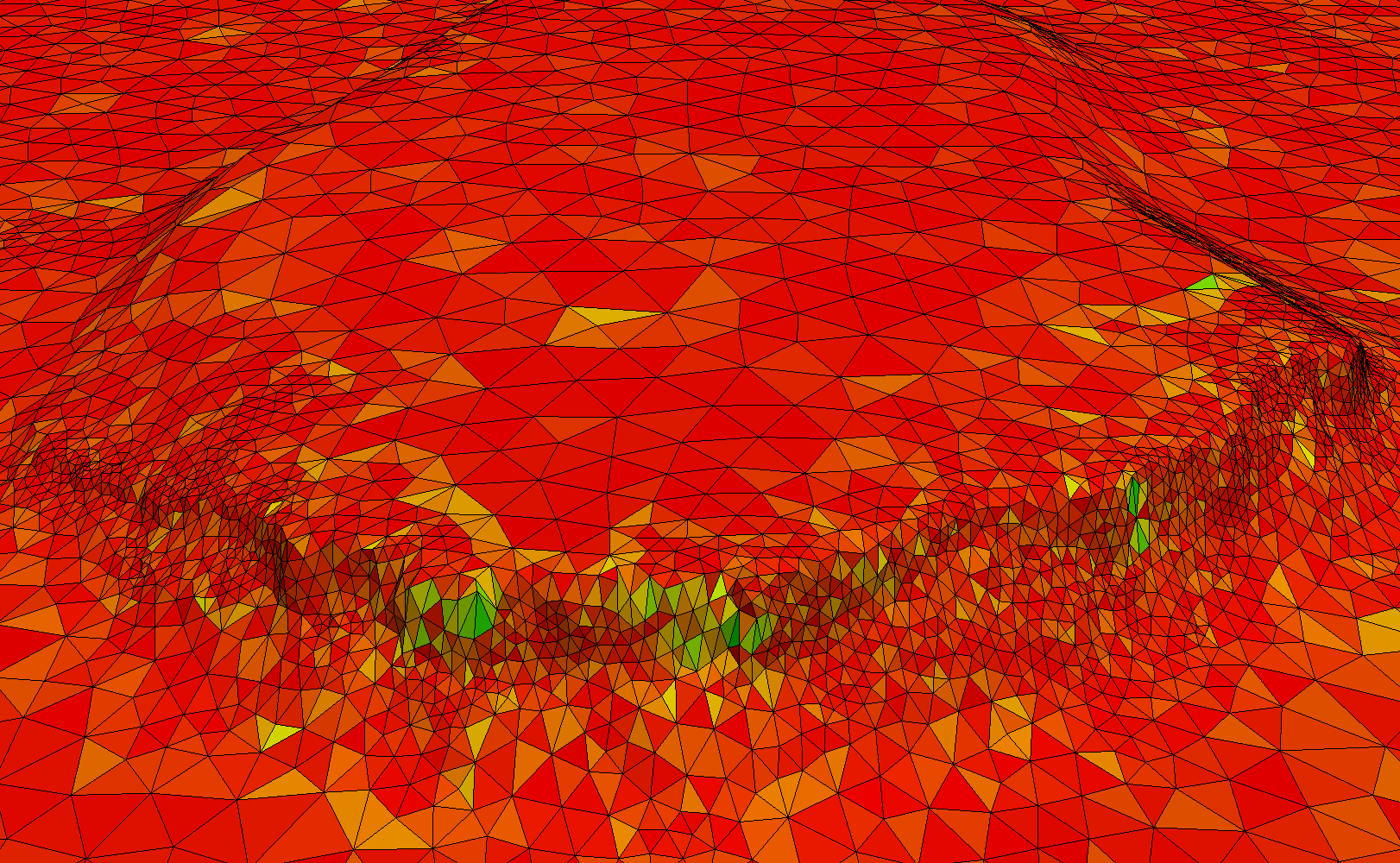}  	 } 	
	\end{tabular}
	\\
	\includegraphics[width=0.2\textwidth]{./qualBarParaview_color} 
	\caption{
		Topography mesh optimization for the Bolund hill: (a) initial distorted surface mesh, and (b) optimized high-quality surface mesh. The surface elements elements are colored according to their quality,  Eq. \eqref{eq:quality}.}
	\label{fig:bolundSurfaceProcess}
\end{figure*} 

\begin{table} 
	\caption{Shape quality statistics (minimum, maximum, mean and standard deviation) for the surface meshes presented in Figure \ref{fig:bolundSurfaceProcess}.}
	\label{tab:bolundSurfaceProcess}
	\centering
	\begin{tabular}{ccccc}
		\hline
		Mesh  & {Min.Q.} & {Max.Q.} &{ Mean Q.} & {Std} \\
		\hline
		Fig. \ref{fig:surfaceMesh1} & 0.17    &  1.00 & 0.89  & 0.12 \\
		Fig. \ref{fig:surfaceMesh2} & 0.50    &  1.00 & 0.96  & 0.05 \\
		\hline
	\end{tabular}
\end{table}

In this section, we analyze the mesh quality improvement in the mesh generation proces for the test case presented in Section 
 \ref{sec:volume}, the Bolund hill. In addition, we analyze the effect of the mesh optimization in the convergence of the solver for different scenarios.

The optimization of the adapted surface mesh
for the Bolund hill  is shown in 
Figure \ref{fig:bolundSurfaceProcess}. In addition, in Table \ref{tab:bolundSurfaceProcess} we illustrate the mesh quality statistics resulting from the detailed procedure.  Figure \ref{fig:surfaceMesh1} shows the initial distorted surface mesh, which has a minimum quality of 0.17. 
It can be observed that non-regular lower quality elements are present  in areas with big slopes of the topography.  
Figure \ref{fig:surfaceMesh2} shows the optimized mesh, with a minimum quality that has been increased to 0.50, and where it can be observed that the distorted elements from the initial mesh have become almost regular all over the domain.

We highlight that the quality of the surface mesh is of  major importance for the volume mesh generation. The surface mesh defines the boundary of the volume mesh, and therefore, a low-quality surface mesh will derive in a low-quality (or invalid) volume mesh.
In addition, in Section \ref{sec:convergenceGeometry} we study the convergence of the mesh to the Bolund geometry, depending on different meshing parameters, including the usage or not of surface mesh optimization.
As it will be observed, improving the mesh quality implicitly results in a reduction of the error of the geometry discretization, since improving the element quality provides an improved distribution of the elements on the surface.

\begin{figure}[!t] 
	\centering
	\hspace{-0.55cm}
	
	\begin{tabular}{cc}
		\subfigure[]{\label{fig:prismTet}
			\includegraphics[width=0.4\textwidth]
			{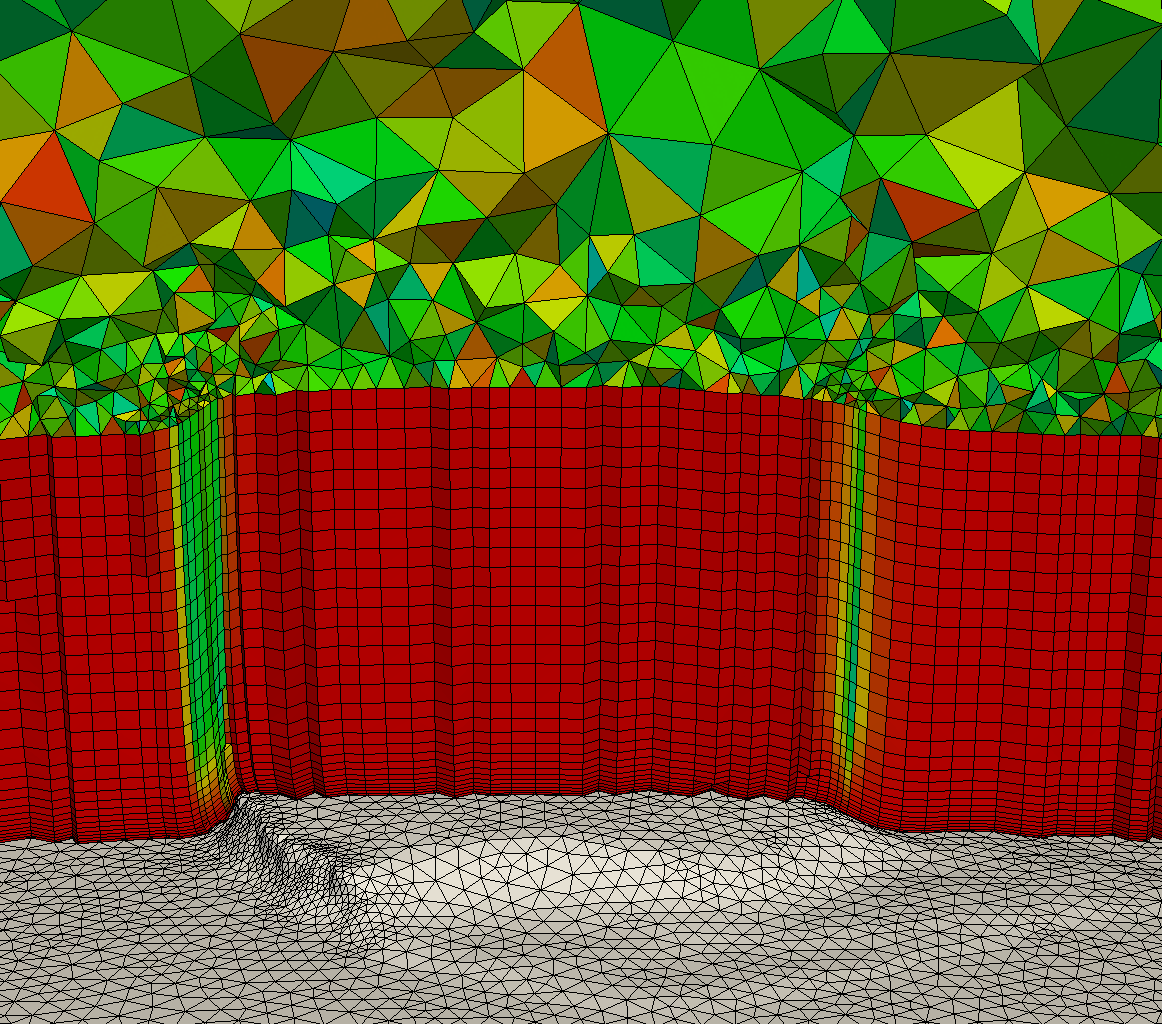}  	 } 
		&
		\subfigure[]{\label{fig:prismTetOpti}
			\includegraphics[width=0.4\textwidth]
			{./bolund_volume_process.0002}  	 } 
	\end{tabular}
	\includegraphics[width=0.2\textwidth]{./qualBarParaview_color} 
	\caption{
		Hybrid ABL mesh for the Bolund hill before \subref{fig:prismTet} and after \subref{fig:prismTetOpti} the hybrid mesh optimization.
		The elements are colored according to their quality. 
	}
	\label{fig:bolundWindMesh_optimization}
\end{figure} 

\begin{table} 
	\caption{Shape quality statistics (minimum, maximum, mean and standard deviation) for the volume meshes presented in Figure \ref{fig:bolundWindMesh_optimization}. }
	\label{tab:bolundWindMesh_optimization}
	\centering
	\begin{tabular}{ccccc}
		\hline
		Mesh  & {Min.Q.} & {Max.Q.} &{ Mean Q.} & {Std} \\
		\hline
		Fig. \ref{fig:prismTet} & 0.10    &  1.00 & 0.85  & 0.23   \\
		Fig. \ref{fig:prismTetOpti} & 0.20    &  1.00 & 0.86  & 0.19   \\ 
		\hline
	\end{tabular}
\end{table}%

Figure \ref{fig:bolundWindMesh_optimization}  illustrates the generated hybrid volume mesh 
from the topography surface mesh illustrated in Figure \ref{fig:surfaceMesh2}.
Figure \ref{fig:prismTet} shows the hybrid mesh obtained after meshing the SBL with prisms and the rest of the domain with tetrahedra.
In Figure \ref{fig:prismTetOpti} we illustrate the final hybrid mesh obtained after the performed quality optimization.
Table \ref{tab:bolundWindMesh_optimization} details the quality of the meshes presented in Figure \ref{fig:bolundWindMesh_optimization}.
The minimum quality of the mesh in Figure \ref{fig:prismTet}  is 0.1, whereas in Figure \ref{fig:prismTetOpti} the minimum quality has increased up to 0.2.
In addition, if no volumetric optimization is performed during the extruding procedure (see Line \ref{alg:optimizePrisms} of Alg. \ref{alg:generatePrismMesh}), the mesh becomes invalid during the sweeping process and the mesh can not be generated.
This issue stresses the importance of the quality optimization framework, specially in regions where the topography features high gradients.

\begin{figure}[!t] 
	\centering
	
	\subfigure[]{\label{fig:alaiz}
		\includegraphics[width=0.7\textwidth]{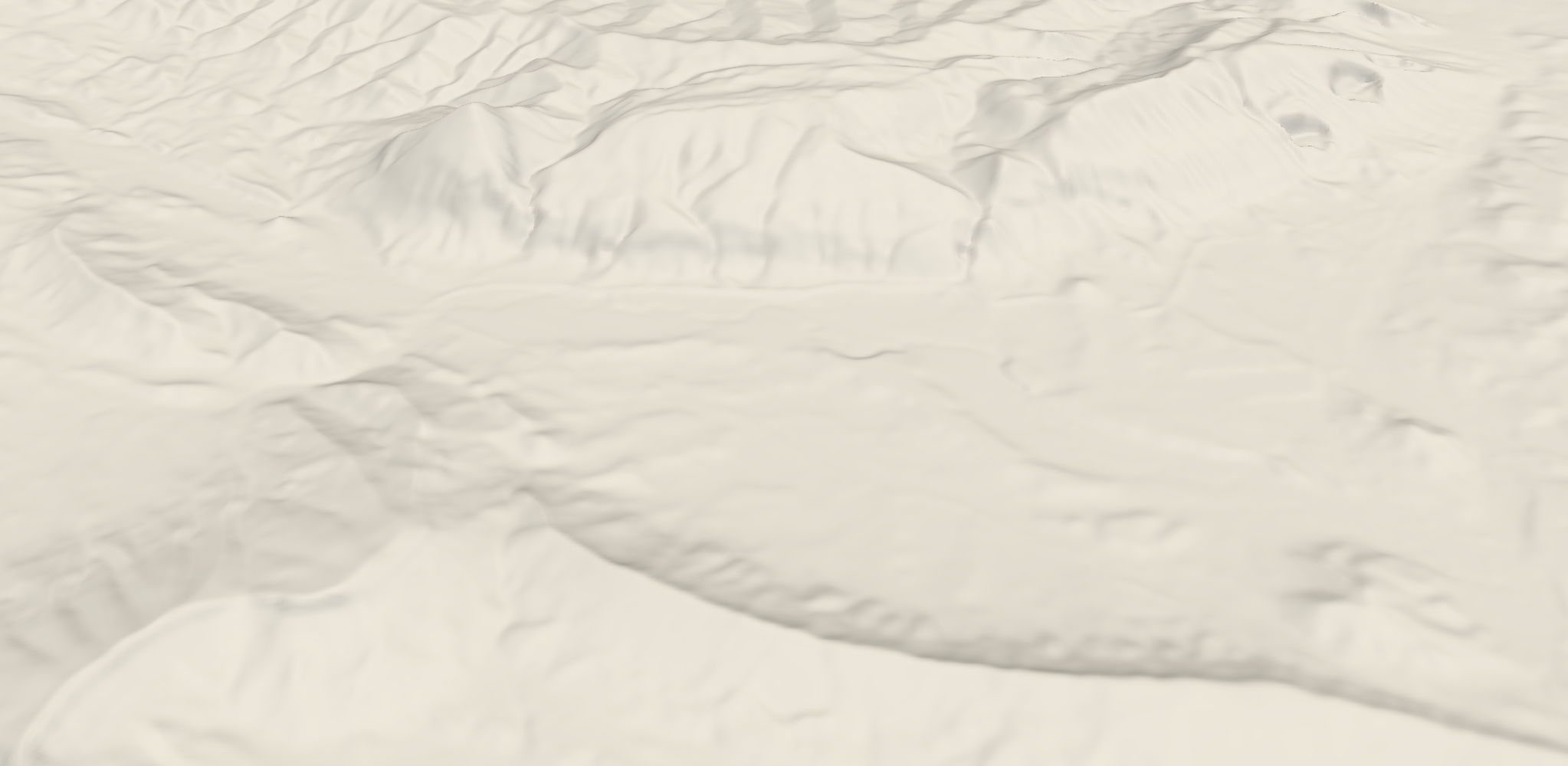}  	 
	}
	\\
	\subfigure[]{\label{fig:pueblaTopo}
		\includegraphics[width=0.7\textwidth]{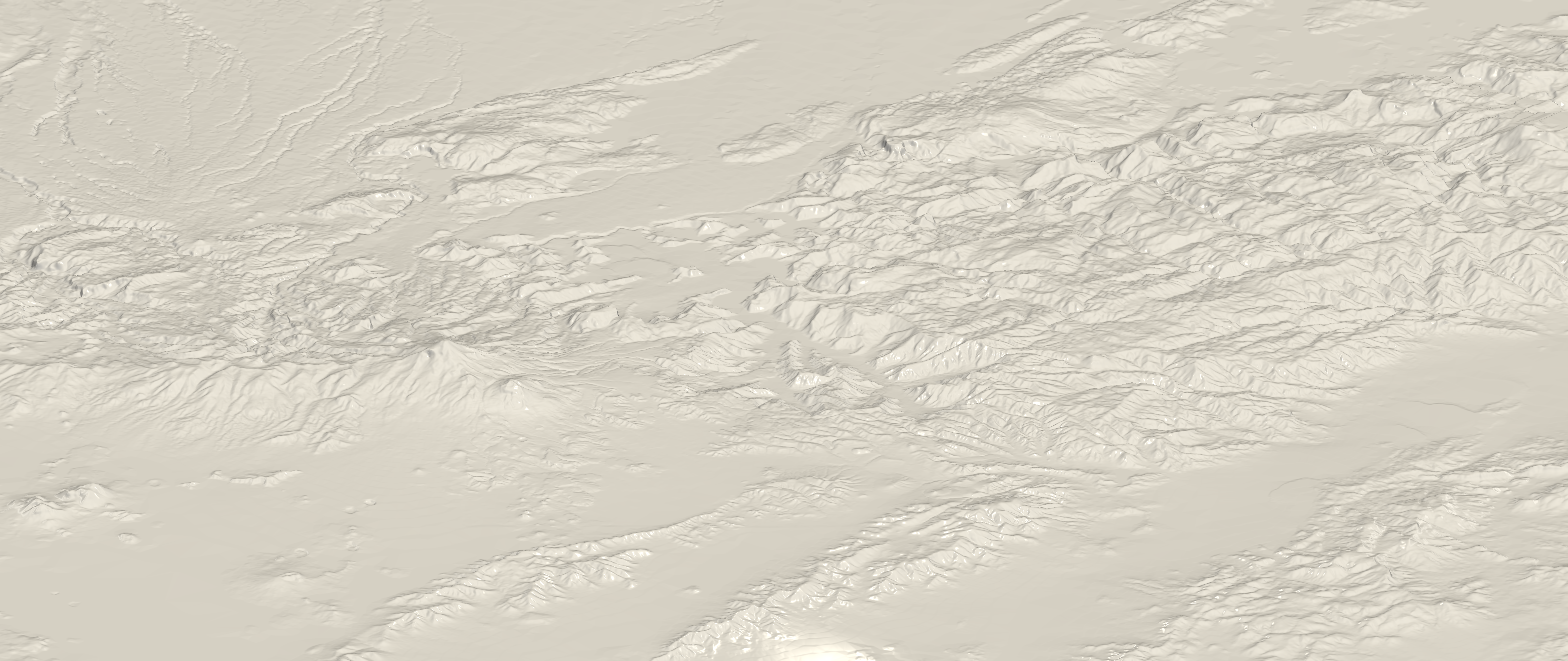}  	 
	}	
	\caption{
		Topography of the:
		\subref{fig:alaiz}	Alaiz mountain range  (Spain), and
		\subref{fig:pueblaTopo}	Puebla topographic scenario (Mexico).
	}
	\label{fig:toposConvergenceOpti}
\end{figure} 

To check the effect of the mesh optimization on the convergence of the solver, we have performed studies in three different topographical landscapes to check the influence of this mesh improvement in the CFD model presented in Section \ref{sec:cfd}.
The influence of the optimization on the resulting meshes can be observed in three different manners.
First, for highly complex scenarios, the simulation is precluded without optimization, due to the lack of quality of the elements in the complex features of the geometry.
As previously highlighted, this is the case of the Bolund scenario for small mesh sizes, where the mesh must be optimized to avoid obtaining invalid elements.

\begin{figure}[!t] 
	\centering
		\includegraphics[width=0.6\textwidth]{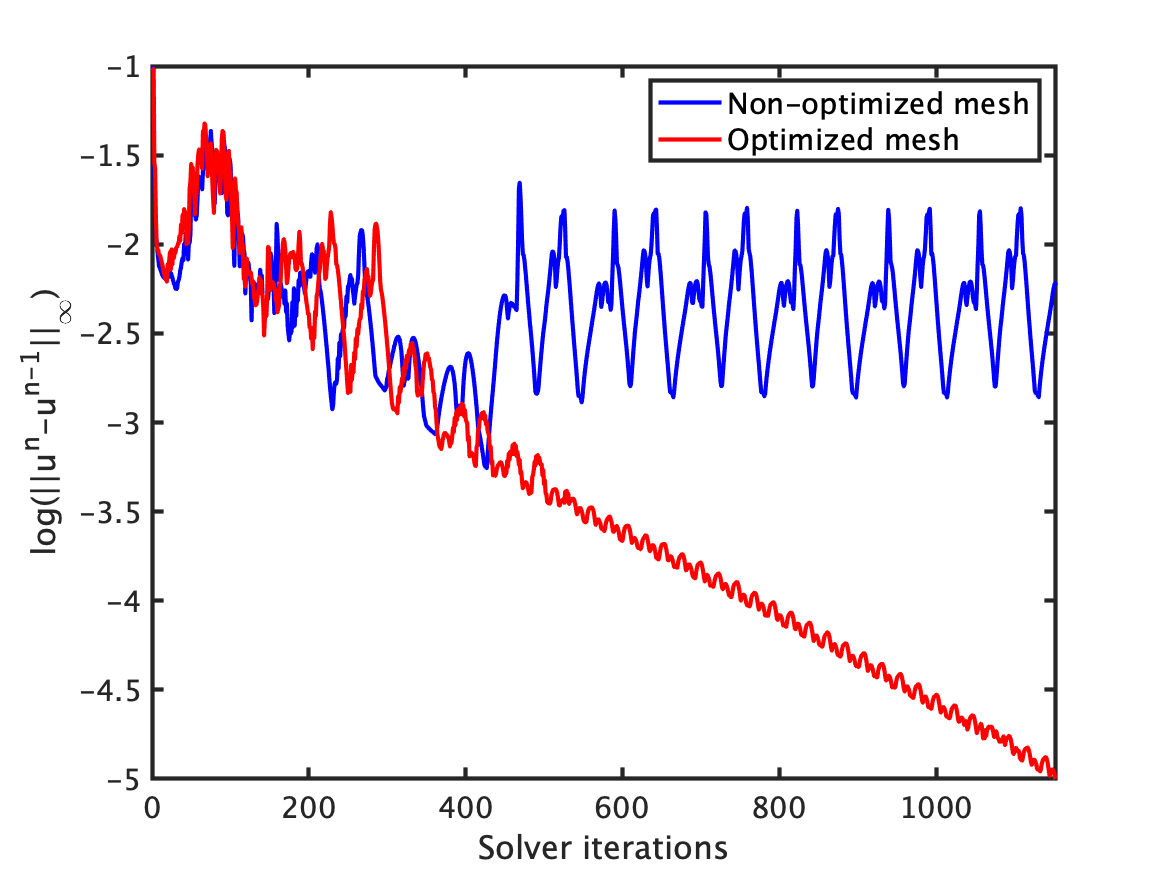}  	 
	\caption{
		Solver convergence on the Alaiz topography, see Fig. \ref{fig:alaiz}, featuring a non-optimized and an optimized mesh.
	}
	\label{fig:alaizConvergence}
\end{figure}

Second, in topography cases where the mesh generated without optimization is valid, two main outcomes can occur.
On the one hand, it can happen that convergence to a steady stat is not achieved. This is the case for the Alaiz mountain range scenario (Spain) presented in Figure \ref{fig:alaiz}. In this scenario, the optimization framework enables the convergence of the solver to a steady state. 
In Fig. \ref{fig:alaizConvergence} we plot for each pseudo-time step, the $\mathcal{L}^{\infty}$-norm of the difference between the current wind velocity with respect to the one from the previous pseudo-time step.
Without optimization, the solution converges to a solver relative tolerance of $10^{-2}$, and then oscillates at this value. 
With the proposed optimization procedure, the solver is able to converge to the desired tolerance (to illustrate here $10^{-5}$), achieving the desired steady state in the RANS simulation.
On the other hand, for some scenarios it can be observed that the optimization allows reaching the steady state  faster than without optimization. That is, the number of pseudo-time steps required to reach the steady state are reduced with optimized meshes, reducing then the required CPU time.
For instance, for the topography illustrated in Figure \ref{fig:pueblaTopo}, a complex scenario located in Puebla (Mexico),
the pseudo-time steps to reach the stationary state are reduced the 50\%
when optimization is enabled.

We highlight that the impact of the mesh optimization depends on the complexity of the terrain, and on the considered mesh size. For simple topographies, the benefit may be a slight reduction of the iterations of the solver. However, for complex topographies it can either enable the simulation (Bolund), enable the convergence of the solver to a steady state (Alaiz), or reduce the number of steps to reach the steady state (Puebla), with the derived reduction of computational cost.
This result is in agreement with the fact that well shaped elements (without small/large angles) improve the matrix conditioning in Finite Elements, as detailed by Babu{\v{s}}ka and Aziz \cite{babuska:angle} and Shewchuk \cite{shewchuk:goodElements}. 


\section{Results}
\label{sec:results}

In this section, we analyze the performance of the presented adaptive hybrid meshing strategy.
First, in Section \ref{sec:exampleBolund} we compare the hybrid meshing strategy without adaptation with respect to the standard semi-structured approach in terms of the degrees of freedom generated by each of the methods.
Second, in Section \ref{sec:convergenceGeometry} we analyze the convergence to the geometry of the proposed mesher with and without adaptation, obtaining quadratic mesh convergence for both cases, and reducing the error when adaptation is enabled.
Third, in Section \ref{sec:convergenceSolution} we analyze the benefits of the proposed approach for ABL simulation using the CFD model presented in Section \ref{sec:cfd}.  We compare it with respect to our semi-structured procedure, attaining quadratic convergence for both scenarios, but resulting in a significant reduction of the required degrees of freedom with the hybrid approach.
Finally, in Section \ref{sec:exampleBadaia} we illustrate the applicability of the approach to simulate the ABL flow over the complex topography from
 the Badaia wind farm in Spain.

Following, we detail the computational time to generate the main meshes presented in this section, all generated on a MacBook Pro with one dual-core Intel Core i7 CPU,  a clock frequency of 3.0 GHz, and a total memory of 16 GBytes.
In Section \ref{sec:exampleBolund}, the computational time to generate the semi-structured mesh in Figure \ref{fig:quaBolundVolume} is 9 seconds, and the time to generate the hybrid mesh in Figure \ref{fig:triBolundVolume} is 17s.
In Sections \ref{sec:convergenceGeometry} and \ref{sec:convergenceSolution} a total of 28 meshes are generated to assess the mesh convergence, 
which feature a small element count,
 and each generated in a few seconds.
Finally, the computational time to generate the hybrid mesh in Figure \ref{fig:badaiaVolume} from Section \ref{sec:exampleBadaia} is 246 seconds.

\subsection{Comparison of structured vs unstructured approaches} 
\label{sec:exampleBolund}

The objective of this work is to propose a new mesher to simulate ABL flows featuring Coriolis effects in complex topographies.
In this context, the motivation of this section is to present a simple example to illustrate the advantages of using the proposed approach against a standard semi-structured ABL mesher in terms of the resulting element count and the flexibility to attain the desired element sizes in the different regions of the  domain.

\begin{figure}[t!] 
	\centering
	\begin{tabular}{cc}
		\subfigure[]{\label{fig:quaBolundVolume}
			\includegraphics[width=0.45\textwidth]
			{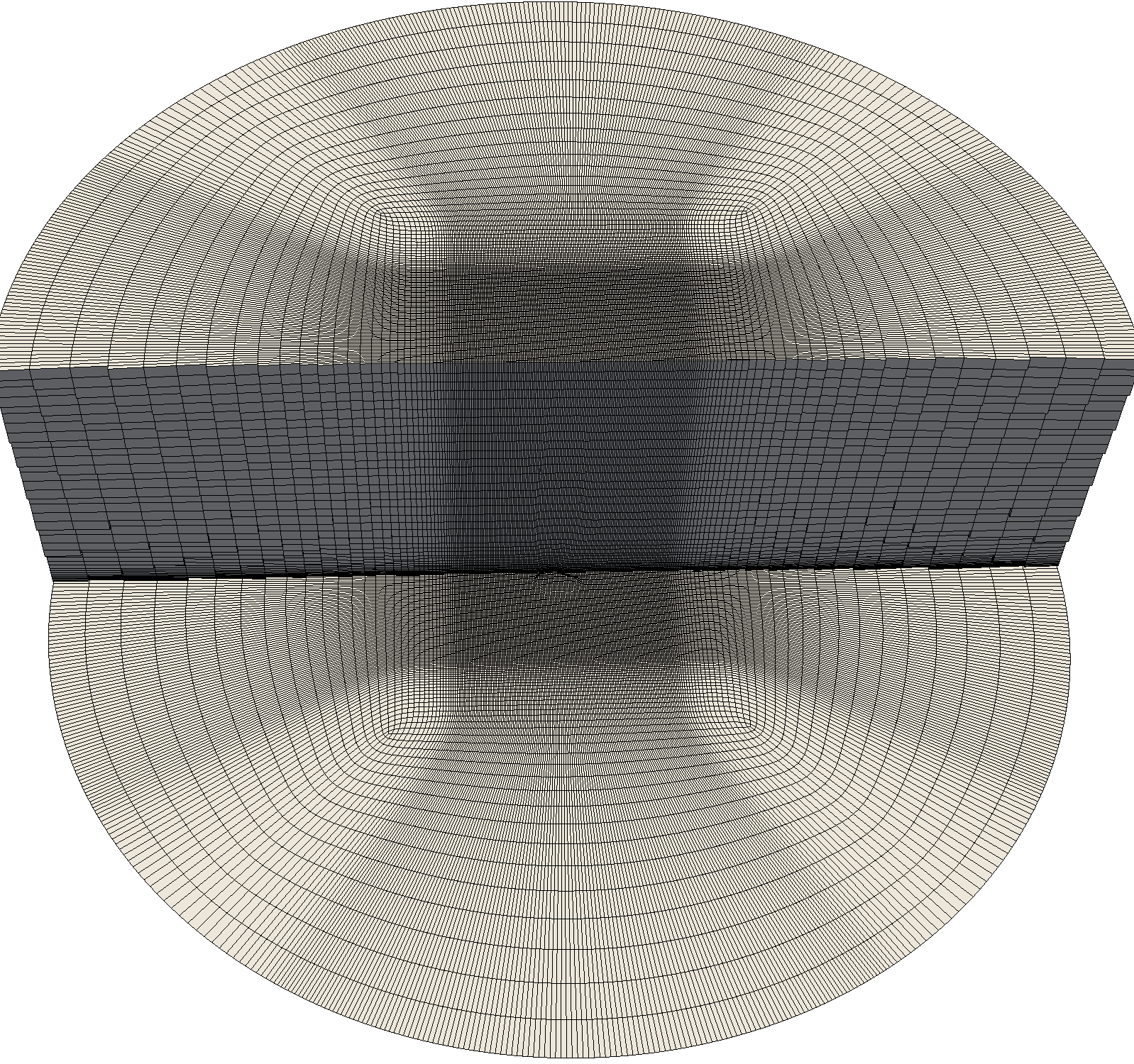}  	 }
		&
		\subfigure[]{\label{fig:triBolundVolume}
			\includegraphics[width=0.45\textwidth]
			{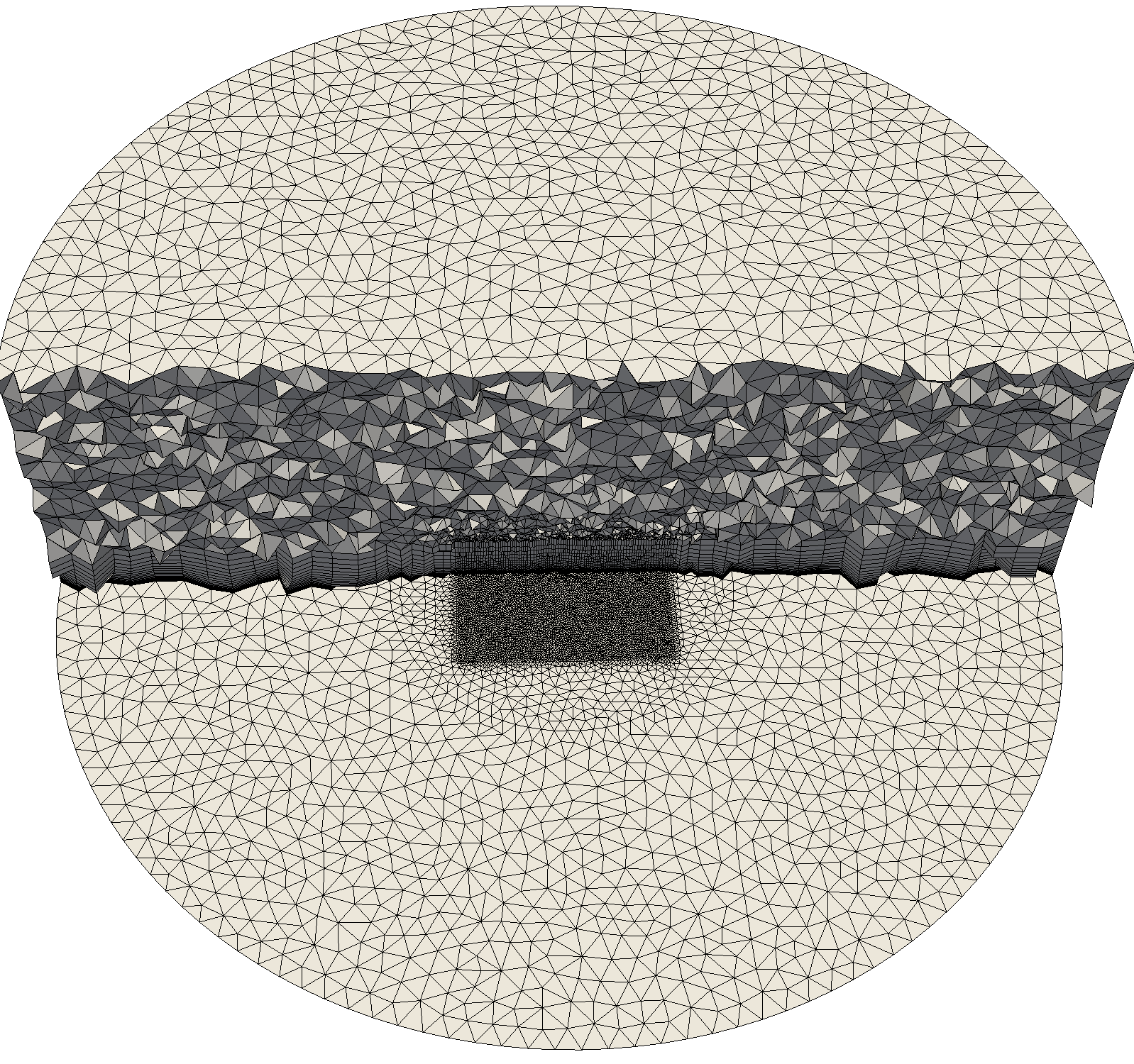}  	 }
	\end{tabular}
	\caption{
		Meshes generated on Bolund hill: \subref{fig:quaBolundVolume} standard semi-structured mesh and
		\subref{fig:triBolundVolume}  hybrid mesh.
	}
	\label{fig:bolundVolume}
\end{figure}

One of the targets of the hybrid mesher is to reduce the degrees of freedom required to discretize a complex topography and the Atmospheric Boundary Layer in simulations featuring Coriolis effects.
To perform this comparison, we use our in-house semi-structured code \xWindMesh\ \cite{gargallo:meshForABLandWindFarms,gargallo2018JCP:WindFarms}, that generates a semi-structured cross-type mesh, see Figure \ref{fig:quaBolundVolume}, which constitutes the standard in
the field of ABL flow simulation \cite{sorensen:hypgrid,sorensen:ellypsys3Dvalidation,openfoam:web,gargallo:meshForABLandWindFarms,avila:CFDframworkWindFarms,zephy,gargallo2018JCP:WindFarms}. 
As previously highlighted, these semi-structured meshers are devoted to generate an hexahedral mesh aligned with the inflow direction.
Thus, a structured quadrilateral mesh is generated in the interest regions, and then an elliptical domain is meshed using a semi-structured approach increasing the mesh size towards the exterior of the domain. Next, this surface mesh is extruded in the vertical direction to discretize the boundary layer, see \cite{gargallo:meshForABLandWindFarms,gargallo2018JCP:WindFarms} for  further details.

Figure \ref{fig:bolundVolume} illustrates the obtained meshes. Figure \ref{fig:quaBolundVolume} shows the generated semi-structured mesh, while Figure \ref{fig:triBolundVolume} we illustrate the hybrid mesh generated with exactly the same parameters than the semi-structured one. That is, the imposed mesh sizing is  the same. 
In this particular test case, featuring a small domain with a small topographic scenario around Bolund hill, the input sizes are: 5 meters in the farm region and 25 in the buffer region. Regarding the vertical discretization, the initial cell size is 0.5 and the geometrical growing ratio is 1.15.
We highlight that in this example, we do not adapt the hybrid mesh to the topography so that both meshers are comparable.

Regarding the number of elements and nodes of the generated meshes, the hexahedral mesh in Figure \ref{fig:quaBolundVolume} is composed by 1.41 Million (M) nodes  and  1.36M elements. 
In contrast, the meshes in Figure \ref{fig:triBolundVolume} is composed by 0.27M nodes and 0.77M elements, from which 0.41M  are prisms and 0.36M are tetrahedra.
In particular, we obtain a mesh with half of the elements and a fifth of the number of nodes but attaining the same required resolution than in the structured approach. 
Recall that the computational cost of a given simulation depends on the number of nodes, since the nodes determine the number of unknowns to solve.
Thus, a fifth of the number of nodes is a significant reduction of the computational cost and of the size of the matrices involved in the solvers, which have a dimension of the order of the number of nodes.
However, we have the additional advantage of being able to adapt the mesh to ensure that the element size on the surface is the desired and to capture the curvature of the geometry. 
The effect of this adaptation will be observed in Section \ref{sec:convergenceSolution}.

\subsection{Mesh convergence to geometry: uniform vs adaptive refinement}
\label{sec:convergenceGeometry}

\begin{figure}[!t] 
	\begin{center}
		\includegraphics[width=0.475\textwidth]
		{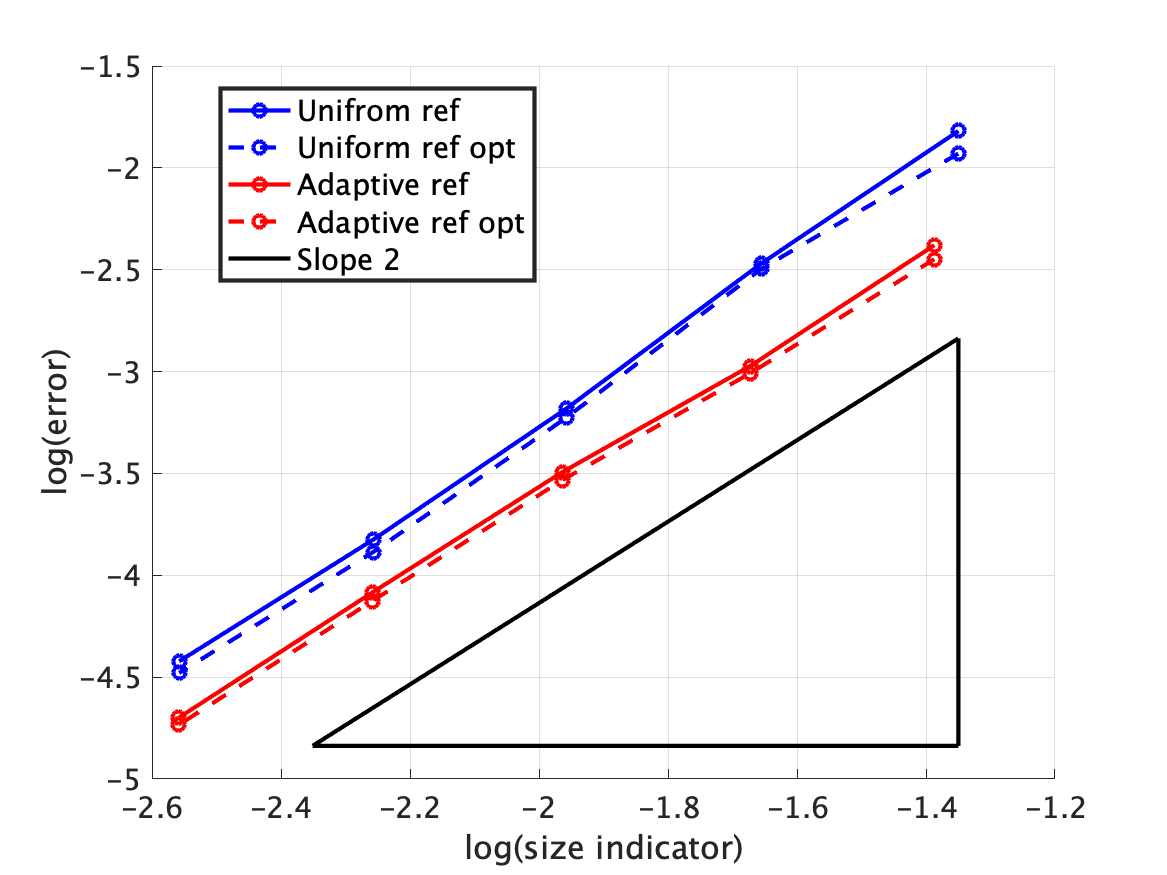}  	
	\end{center}
	\caption{Logarithm of the geometry error with respect to the logarithm of the size indicator for the uniform and adaptive refinement. 
	}
	\label{fig:uniVsAdaptConv}
\end{figure}

\begin{figure}[!t] 
	\begin{center}
		\includegraphics[width=0.475\textwidth]
		{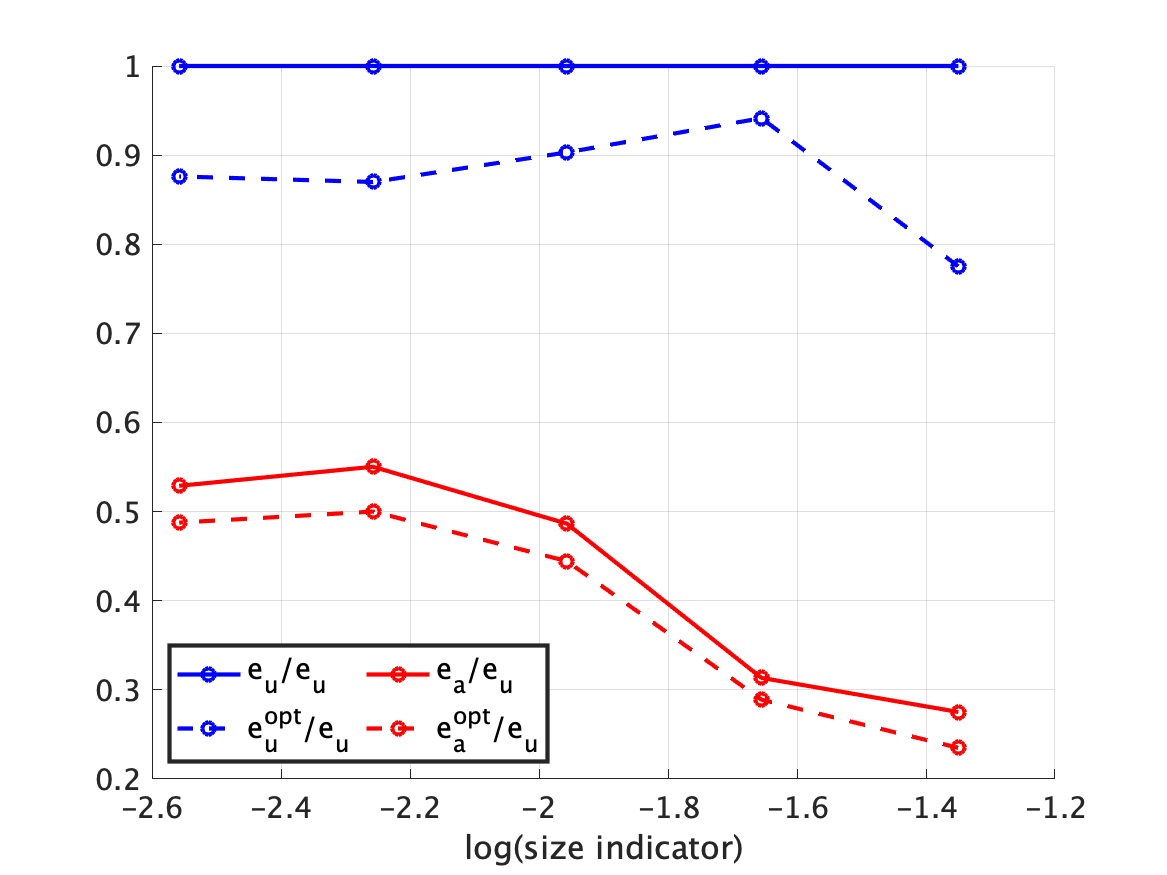}  	
	\end{center}
	\caption{Error reduction by enabling mesh optimization and adaptivity: logarithm of the size indicator versus the ratio of mesh error with respect to the uniform mesh error. Notation: $e_u$ is the $\mathcal{L}^2$ error of the uniform mesh, $e_u^{opt}$ is the error of the optimized uniform mesh, $e_a$ of the adapted mesh, and $e_a^{opt}$ of the optimized adapted mesh.}
	\label{fig:uniVsAdaptConv_red}
\end{figure}

In this section, we focus  on the scenario of the Bolund peninsula to study the influence on the geometric accuracy of the proposed adaptive procedure in contrast to uniform mesh refinement. 
The advantages in terms of degrees of freedom of unstructured triangle surface meshes (hybrid volume meshes), against semi-structured quadrilateral surface meshes (hexahedral volume meshes), have already been analyzed in Section \ref{sec:exampleBolund}. 
Thus, in this section we focus strictly on unstructured triangle meshes, to analyze the advantages delivered by the adaptation of the triangle mesh to the topography.

On the one hand, we  generate a sequence of triangle topography meshes starting by a uniform mesh size of 40 m, and dividing by 2 the mesh size. In particular, we generate 4 meshes of edge length approximately 40, 20, 10, 5 and 2.5 meters. The topography mesh is generated using Algorithm \ref{alg:adaptiveSurface}, with $h_{max}=h_{min}$.
On the other hand, we  also generate the meshes with the same maximum size, but allowing the mesh to be refined according to the two metrics up to 1/4 of the maximum mesh size. That is, $h_{max}=40,20, 10, 5, 2.5$ and $h_{min}=10, 5, 2.5,1.25,0.625$.

For each generated mesh, we also generate a new one enabling the surface optimization procedure, see Section \ref{sec:optimizationSurface}, to illustrate the benefits of an improved mesh quality not only for the solver, but also to improve the surface discretization.
All the generated meshes have approximately around 0.2 of minimum quality before the optimization procedure, which is increased above 0.5 once the mesh is optimized.
The optimization only takes into account the shape of the elements, not the accuracy of the representation of the topography. However, big slopes on the terrain produce low quality configurations and thus, optimizing the mesh quality results in practice in an improvement of the accuracy of the discretization of the topography.

As detailed in the refinement analysis performed in \cite{arnold2000locally}, herein we use $\xNumNodes^\frac{-1}{d}$ as an indicator of the mesh size for locally refined meshes, being $\xNumNodes$ the number of nodes of the mesh and $d$ its  dimension ($d=2$ for triangles, $d=3$ for tetrahedra).
For each mesh we compute the $\mathcal{L}^2$ error of the surface meshes to the given topography. 
In Figure \ref{fig:uniVsAdaptConv} shows a  plot of the logarithm of the error against  the logarithm of the size indicator.
 We show the plot for both uniform and adaptive refinement, for both non-optimized and optimized configurations.
As it can be observed, the theoretical quadratic convergence rates are obtained with all the approaches.
However, the error with the adaptive procedures is reduced with respect to the uniform ones. 
It is also of interest to analyze the indirect benefit produced by the optimization. For both the uniform and adaptive meshes, the optimization still maintains the quadratic convergence and, in addition, reduces the geometry discretization error. 

Figure \ref{fig:uniVsAdaptConv_red} plots the reduction of error between the uniform and adaptive cases enabling or not optimization with respect to the uniform mesh as reference.
For each approach, we plot its error divided by the error of the uniform procedure without optimization, illustrating the percentage of reduction of error with respect to the  non-optimized uniform mesh refinement, which constitutes the standard approach in the literature.
The dashed red line corresponds to the improved procedure proposed in this work,  featuring both local adaptivity and surface mesh optimization.
We can observe that in the analyzed case the error  is reduced almost to the $50\%$ when using the adaptive refinement, and is always reduced over the $50\%$ when both adaptivity and optimization are used. 
Thus,  adapting  and optimizing the mesh results in a more efficient use of the degrees of freedom in terms of the accuracy of the topography discretization.

\subsection{Mesh convergence to solution: semi-structured vs hybrid adapted strategies}
\label{sec:convergenceSolution}

In this section, we analyze the error and the convergence to the solution to the ABL flow model in Section \ref{sec:cfd} of the standard semi-structured ABL meshers with respect to the proposed hybrid approach.
Similarly to Section \ref{sec:exampleBolund}, to perform this comparison we use our in-house semi-structured code \xWindMesh\ \cite{gargallo:meshForABLandWindFarms,gargallo2018JCP:WindFarms}.
The comparison is performed in the Bolund hill, and we measure the mesh convergence to the velocity in a sequence of refined meshes. 

We perform the study with the same  mesh sizing parameters for both approaches, except in the minimum allowed mesh size for the hybrid approach.
Since the objective is to study the influence of the adaptation of the mesh to the topography, the vertical discretization will be fixed for all the tests. 
Regarding the vertical discretization, the initial cell size is 1m and the growing ratio is 1.1.
In particular, 3 meshes will be generated for each approach changing the discretization of the topography,  featuring edge lengths of 20, 10 and 5 meters.
For the adaptive approach, the minimum edge length will be allowed to locally refine the mesh up to one half of the maximum length. 
The error will be measured with respect to a fine mesh of 1m resolution.

\begin{figure*}[!t]
	\centering
	\subfigure[]{
		\label{fig:bolund_veloc}
		\includegraphics[width=0.3\textwidth]{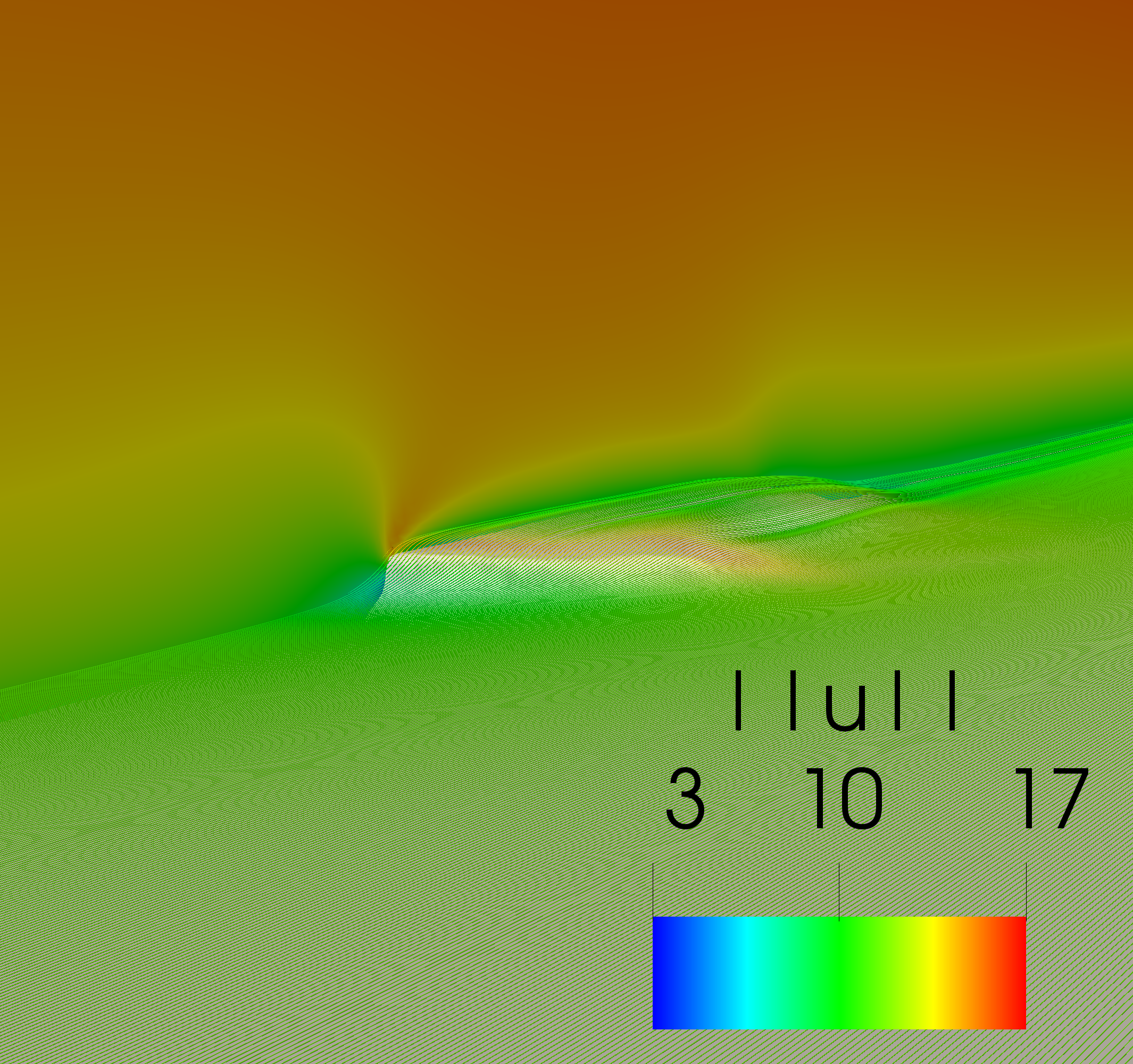} 
	}
	\begin{tabular}{cc}
		\subfigure[]{
			\label{fig:bolund_conv_veloc}
			\includegraphics[width=0.475\textwidth]{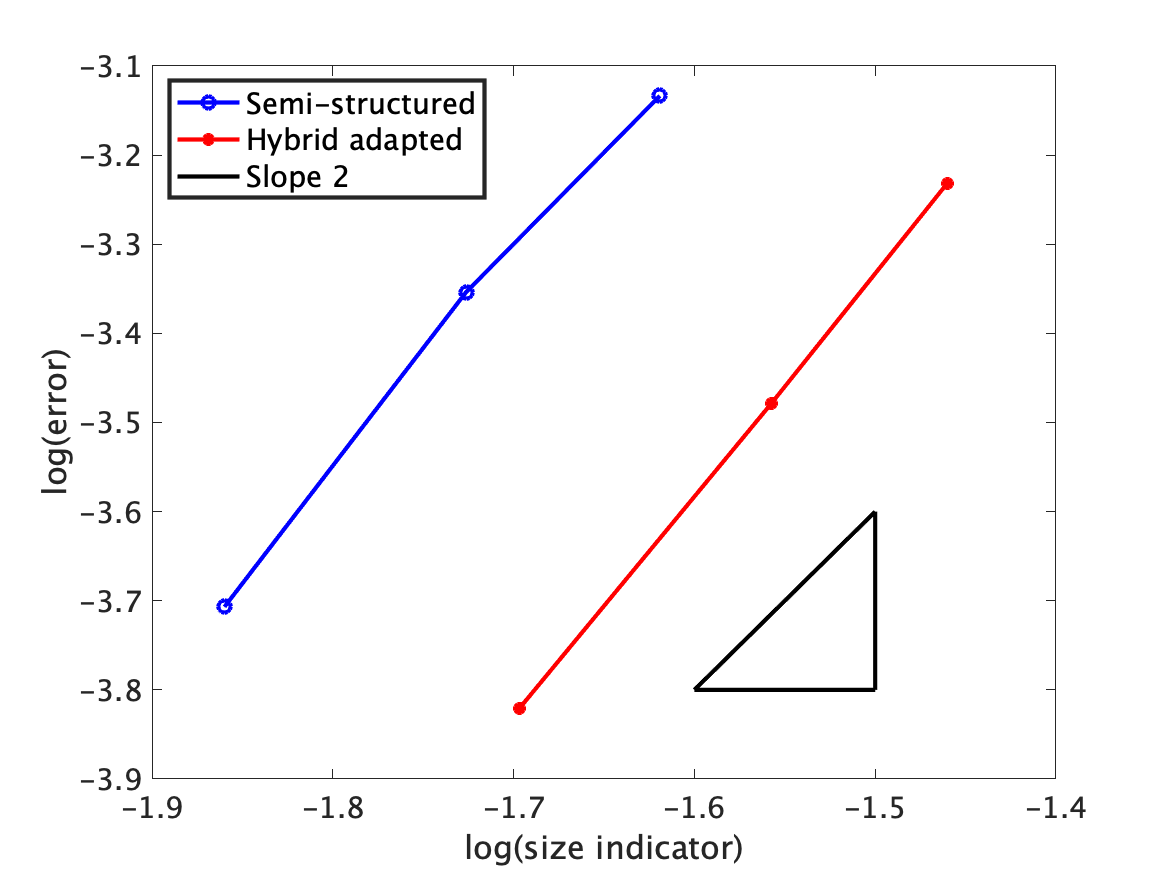} 
		}
		&\hspace{-0.5cm}
		\subfigure[]{
			\label{fig:bolund_conv_veloc_3d}
			\includegraphics[width=0.475\textwidth]{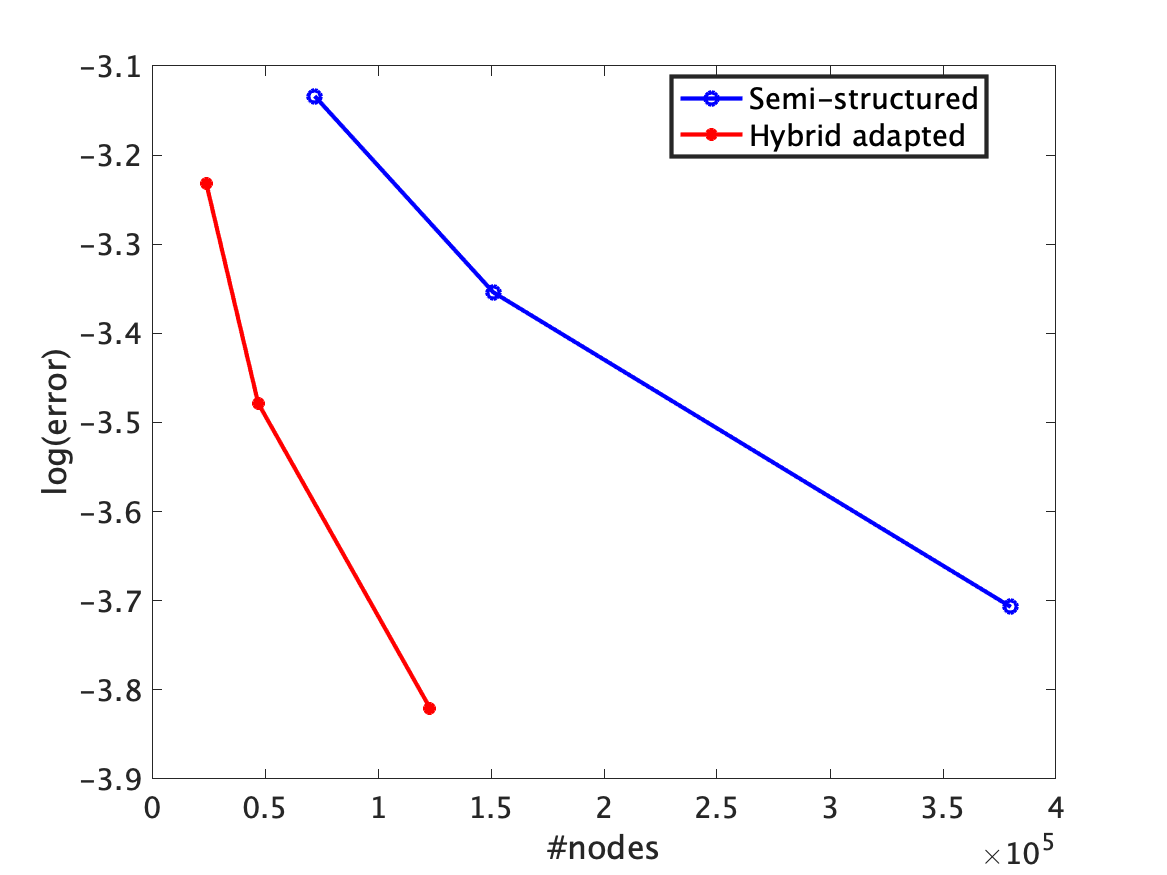} 
		}
	\end{tabular} %
	\caption{
		\subref{fig:bolund_veloc} Reference wind speed over the Bolund hill.
		\subref{fig:bolund_conv_veloc} Logarithm of the velocity error with respect to the logarithm of the size indicator for the semi-structured and  hybrid approaches. 
		\subref{fig:bolund_conv_veloc_3d} Logarithm of the velocity error with respect to the number of nodes of the volume mesh.
	}
	\label{fig:bolund_u_conv}
\end{figure*}

\begin{table} 
	\caption{Statistics of the mesh convergence to solution from Figure \ref{fig:bolund_u_conv}. The ratio of error (nodes) is computed dividing the error ($\#$nodes) of the hybrid with respect to the semi-structured approach.}
	\label{tab:bolundConv}
	\centering
	\begin{tabular}{cccccc}
		\hline
		$h_{max}$  &Type& $\#$nodes & Rel. Error & Ratio error & Ratio nodes \\
		\hline
		\multirow{2}{*}{20} 
		&Hybrid &$24\cdot10^{3}$ &  $0.586\cdot10^{-3}$ &\multirow{2}{*}{0.770}&\multirow{2}{*}{0.33}\\
		& Structured & $72\cdot10^{3}$& $0.735\cdot10^{-3}$ &  &\\ 
		\hline
		\multirow{2}{*}{10} 
		&Hybrid &$47\cdot10^{3}$ &  $0.332\cdot10^{-3}$ &
		\multirow{2}{*}{0.751}&\multirow{2}{*}{0.31}\\
		& Structured & $151\cdot10^{3}$& $0.442\cdot10^{-3}$&  &\\ 
		\hline
		\multirow{2}{*}{5} 
		&Hybrid &$123\cdot10^{3}$ &  $0.151\cdot10^{-3}$ &\multirow{2}{*}{0.799}&\multirow{2}{*}{0.32}\\
		& Structured & $380\cdot10^{3}$& $0.196\cdot10^{-3}$&	  &\\ 
		\hline
	\end{tabular}
\end{table}

Figure \ref{fig:bolund_veloc} shows a slice along the Bolund hill where the wind speed is illustrated.
All the computations have been run in MareNostrum4 \cite{MN4}. 
Figure \ref{fig:bolund_conv_veloc} shows the convergence of both the semi-structured and hybrid adaptive approaches, plotting logarithm of the  size indicator with respect to the logarithm of the $\mathcal{L}^2$ norm of the error of the velocity.
We can observe that both approaches have asymptotic quadratic convergence.
However, the hybrid adaptive procedure reduces the error with respect to the semi-structured approach for the same mesh size. 
Figure \ref{fig:bolund_conv_veloc_3d} plots the logarithm of the error against the number of nodes of the volume meshes, illustrating the reduction of degrees of freedom of the hybrid approach with respect to the semi-structured.
Table \ref{tab:bolundConv} details the mesh and error statistics of the stated cases. 
It can be observed that for each $h_{max}$, the hybrid adaptive approach (allowed to be refined up to $h_{min}=h_{max}/2$), is able to always reduce a $20\%$ the error of the structured approach while only using a $30\%$ of the nodes. 

This study shows the advantages of the hybrid adaptive approach with respect to the standard semi-structured strategies. 
We have to take into account that for each node, 6 unknowns ($\bu\in\xreal{3}$, $p\in\xreal{}$, $k\in\xreal{}$, and $\eps\in\xreal{}$) have to be computed, see Section \ref{sec:cfd}. Thus, requiring up to the 30\% of the nodes (while reducing a 20\% of error), is even more significant once taken into account the problem to be solved.

\subsection{Badaia topography: simulation of ABL flow on a complex scenario}
\label{sec:exampleBadaia}

\begin{figure*} 
	\centering
	\begin{tabular}{cc}
		\subfigure[]{\label{fig:badaia_topo}
			\includegraphics[width=0.4\textwidth]
			{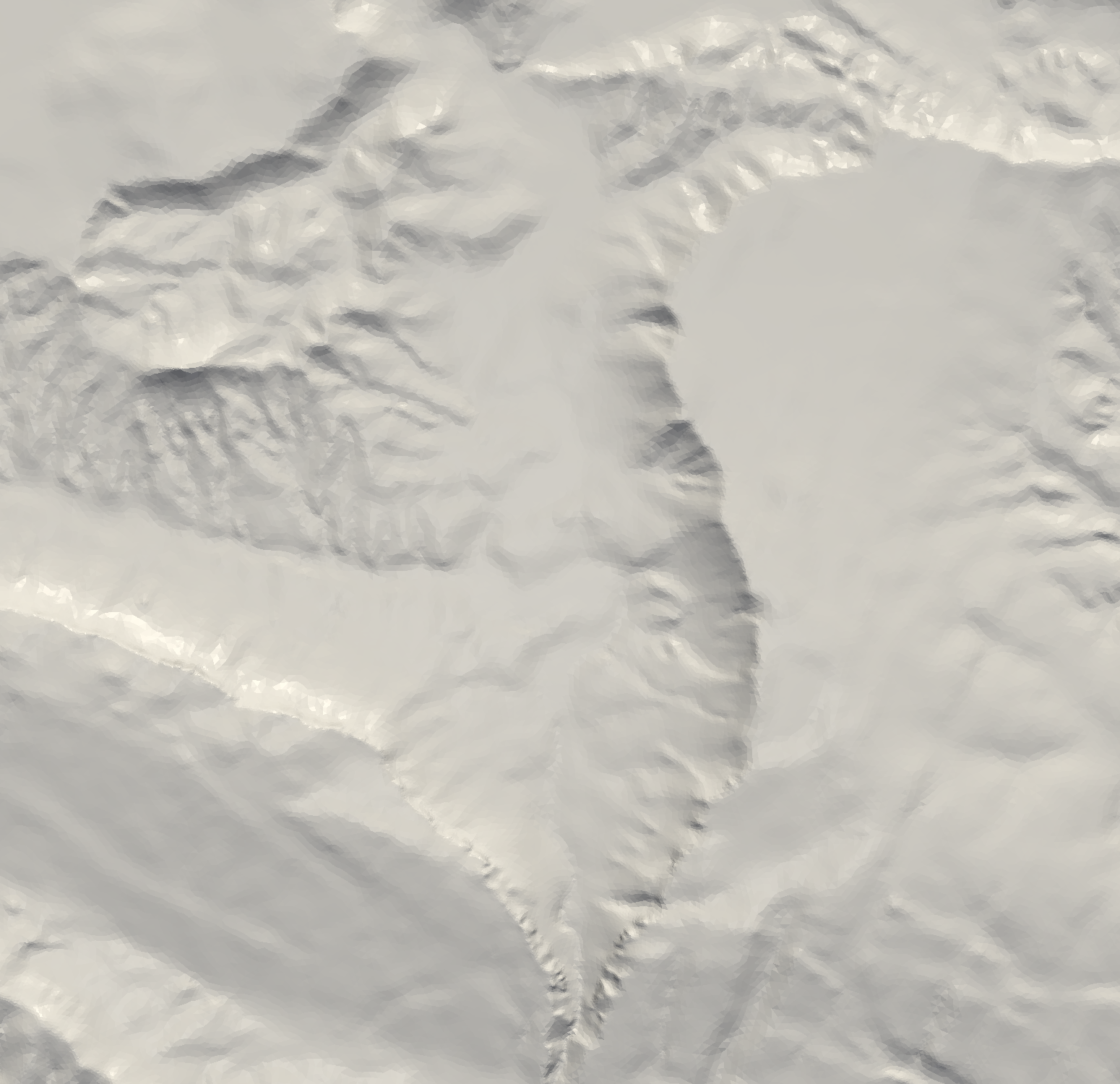}  	 }
		&
		\subfigure[]{\label{fig:badaia_topo_mesh}
			\includegraphics[width=0.4\textwidth]
			{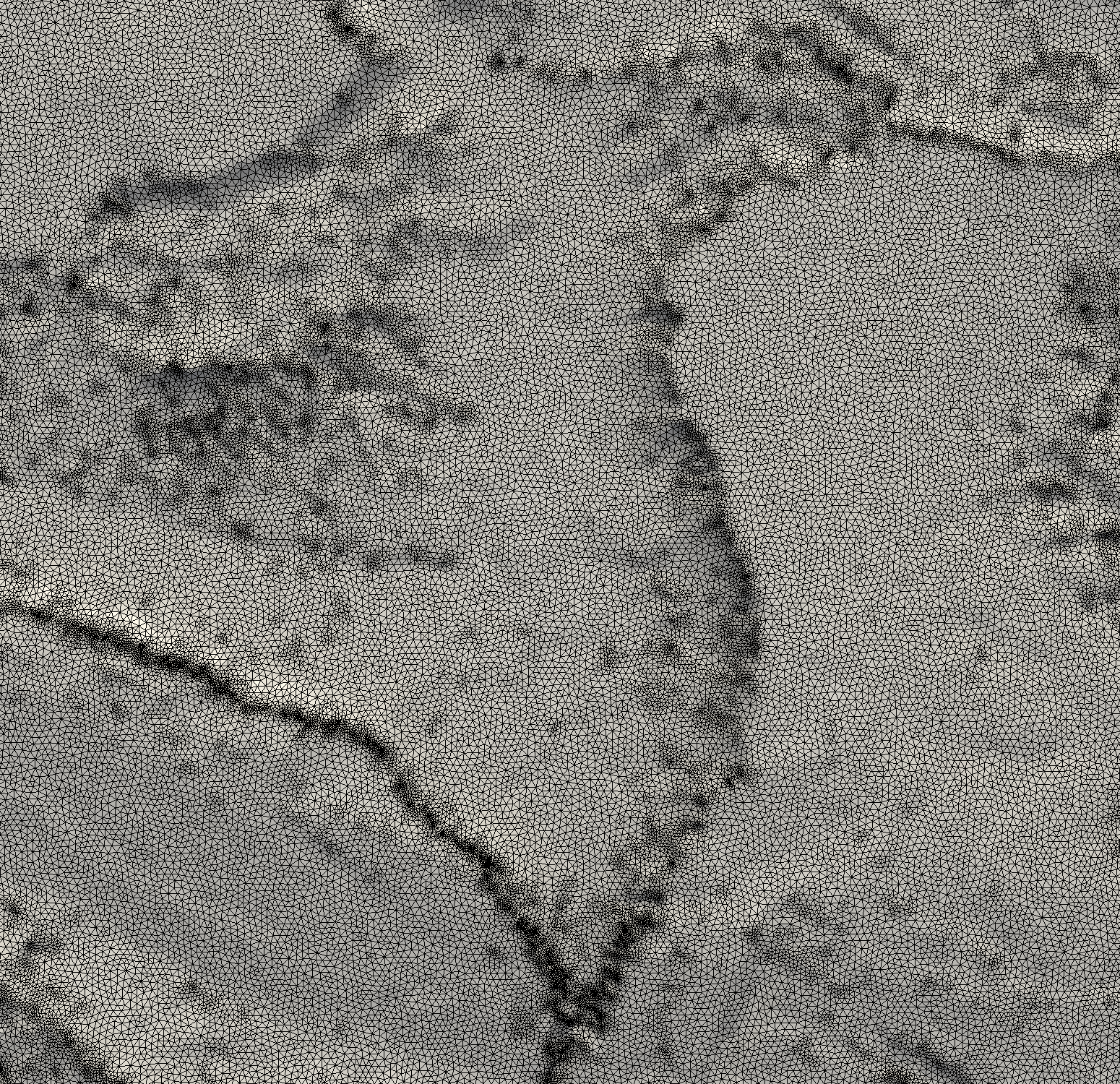}  	 }
	\end{tabular}
	\caption{
		\subref{fig:badaia_topo} Badaia topography, and \subref{fig:badaia_topo_mesh} adapted surface mesh. 
	}
	\label{fig:badaiaSurface}
\end{figure*}

In this section, we present the adapted mesh generated on the Badaia topographic scenario, Figure \ref{fig:badaia_topo}, located in Spain. 
The main objective is to illustrate both the adaptive surface technique and the hybrid volume mesher in a real scenario.
This topography features a valley surrounded by several plateaus and mountains with orographic steps up to 700 meters. 
The minimum height of the topography is 388 meters over the sea level, and the maximum 1098 meters.
The area of the target meshed domain is of 20x20km$^2$, and the top ceiling of the mesh is located at 2km over the highest topographic point.

The surface meshing procedure generates a triangle surface mesh with elements featuring edges of at most 75 meters, 
which are allowed to be refined up to 10 meters to capture the geometry curvature. 
The generated surface mesh is presented in Figure \ref{fig:badaia_topo_mesh}.
The initial mesh is composed by 41531 nodes and 82379 elements. The mesher performs 4  cycles of the refinement procedure presented in Section \ref{sec:surface}. The final surface mesh is composed by  62469 nodes and 124172 elements.
After the adaptive procedure the obtained mesh has a minimum elemental quality of 0.2, which is improved up to 0.24 after performing the optimization procedure presented in Section \ref{sec:optimizationSurface}.

\begin{figure*} 
	\centering
	\subfigure[]{\label{fig:badaia_vol1}
		\includegraphics[width=0.8\textwidth] 
		{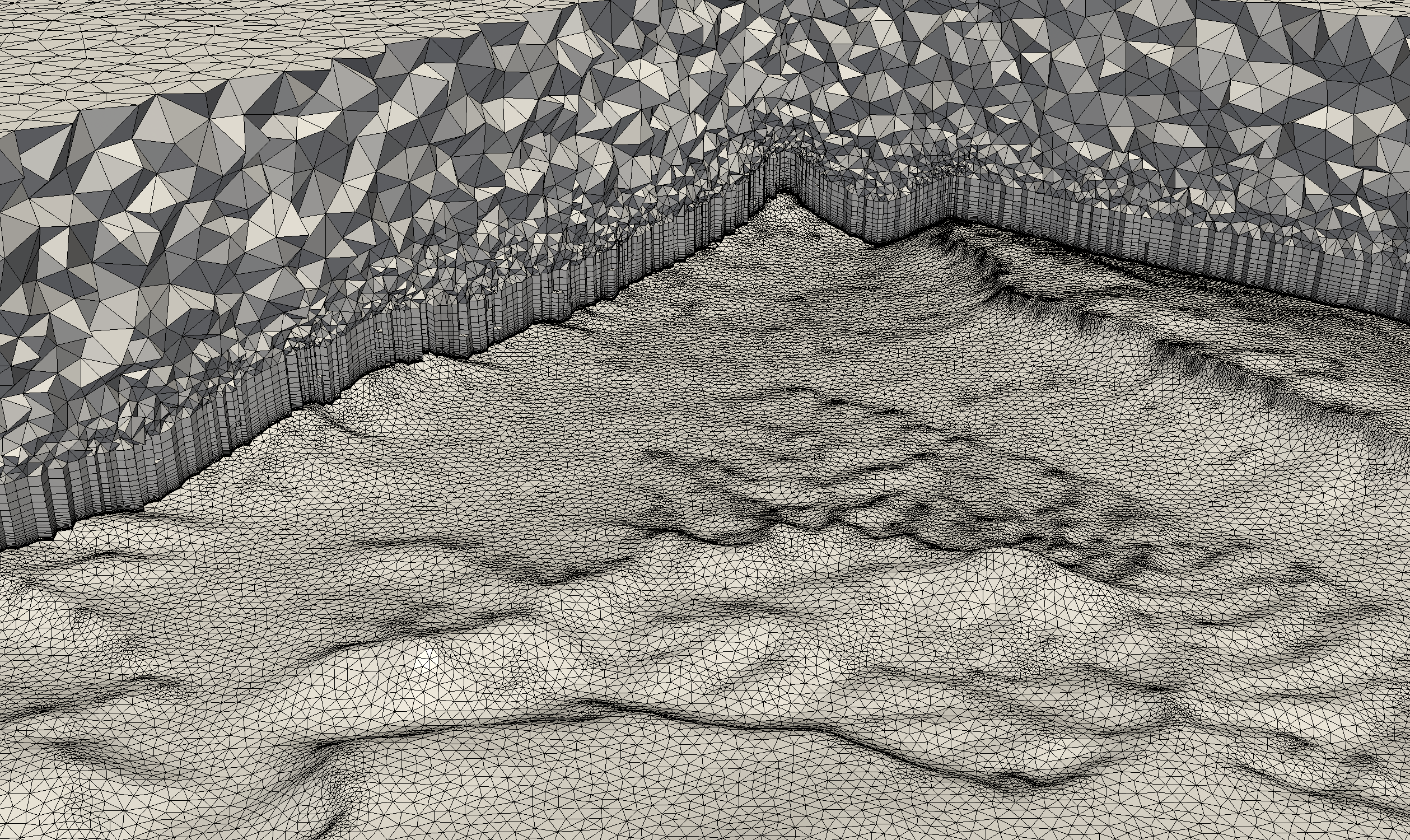}  	 }
	\\
	\subfigure[]{\label{fig:badaiaCFD}
		\includegraphics[width=0.95\textwidth]{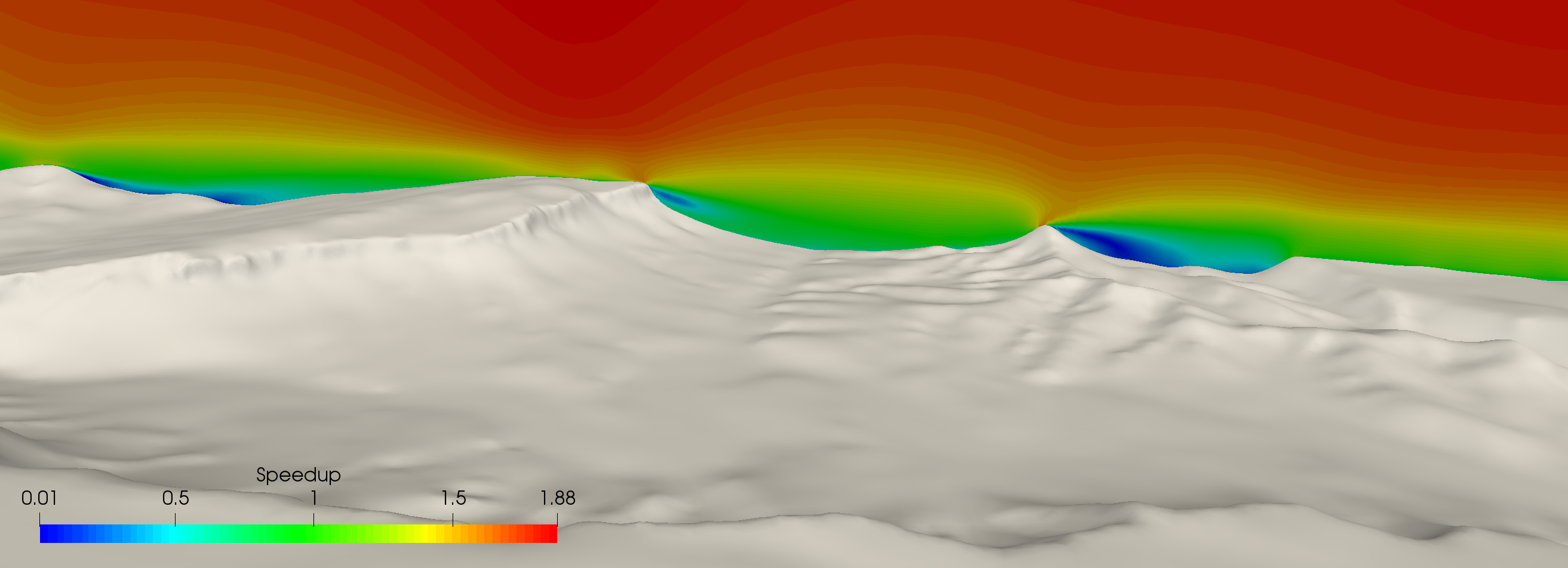}	}
	\caption{
		Badaia scenario:
		\subref{fig:badaia_vol1} ABL adapted hybrid mesh
		and
		\subref{fig:badaiaCFD} velocity speedup with respect to a reference upwind point.
	}
	\label{fig:badaiaVolume}
\end{figure*}

In Figure \ref{fig:badaia_vol1} the volume mesh is presented.
The generated volumetric mesh is composed by 2.2M nodes and 5.6M elements, from which 3.7M are prisms and 1.9M are tetrahedra.  
The mean quality of the generated mesh is 0.84, the minimum is 0.08 and the standard deviation is 0.24.
After the mesh optimization, the minimum mesh quality is improved up to 0.11, the mean remains as 0.84 and the standard deviation is reduced to 0.20. 

In  Figure \ref{fig:badaiaCFD} we illustrate the applicability of the generated meshes for simulation with the model presented in Sec. \ref{sec:cfd}. 
The  simulation has been run in MareNostrum4 \cite{MN4} using 512 cores.
In particular, we illustrate the wind  speedup with respect to a reference point upwind.
The wind inflow direction is set on the left of the domain.
The adaptive approach has locally refined the mesh up to $h_{min}=10$ meters where the topography requires so, whereas the rest of the domain is kept at 75 meters.
Meshing this domain with a semi-structured approach with a resolution of 10 meters to capture the same features than the adaptive approach would require 143M nodes, and 141M elements. This illustrates the potential of the current mesher to simulate large scenarios, and the reduction of the cost of the derived simulation.

\section{Concluding remarks} 
\label{sec:conclusions}

In this work, we present a new hybrid meshing procedure to discretize the Atmospheric Boundary Layer and simulate ABL flows featuring Coriolis force  and complex topographies.
It is a specific purpose mesher to simulate ABL flows, and together with the mesher, the CFD model is detailed, analyzing the mesh requirements derived from the model.
Several contributions in ABL mesh generation have been presented. First, we propose a smooth modeling of the topography that allows queering first and second order derivatives of the geometry for metric computation.
In addition, we incorporate the notion of metric complexity to resolve the curvature metric with the desired number of degrees of freedom.
Second, a new metric-based adaptation procedure to generate the surface mesh of the topography is proposed. This adaptation procedure allows resolving the topography with the desired mesh size and capturing the curvature of the terrain, prescribing a coarser mesh size away from the interest region and locating the nodes to discretize the features of the topography.

Third, we propose to discretize the ABL with a hybrid mesh, featuring prisms in the SBL and tetrahedra in the rest of the domain. Thus, the tensor-structure of prisms is exploited to discretize the boundary layer, and  tetrahedra are exploited to allow a flexible mesh size transition away from the SBL.
In addition, we detail an optimization process for hybrid meshes and analyze the effect of the optimization process on the generated meshes.
It has been observed that for complex topographies the optimization can either enable the simulation, enable the convergence of the solver to a steady state, or reduce the number of steps to reach the steady state, with the derived reduction of computational cost.

To conclude, we present a detailed analysis of the developed meshing strategies. We compare the hybrid adaptive approach with respect to the standard semi-structured meshing strategies, illustrating the advantages in terms of degrees of freedom to converge to both the geometry and to the solution of the ABL flow model. Quadratic convergence to both the solution and the geometry is obtained,  using a 30\% of the degrees of freedom and reducing a 20\% the error with respect to a semi-structured mesh.
In addition, we also apply the proposed approach to mesh and simulate a complex  scenario located in Spain, the topography surrounding the Badaia wind farm.

In the future, regarding the application of the mesher for wind farm resource assessment, 
we would like to combine the capabilities of the proposed mesher to discretize the ABL with complex topographies with the wind farm meshing strategy presented in \cite{gargallo:meshForABLandWindFarms,gargallo2018JCP:WindFarms}.
In particular, the structure of the SBL can be exploited to embed the actuator discs meshes similarly to the approach in  \cite{gargallo:meshForABLandWindFarms,gargallo2018JCP:WindFarms}, resulting in an automatic procedure to mesh onshore wind farms where the hybrid meshing technique is exploited to reduce the required degrees of freedom for the simulation.

\section*{Acknowledgments}
This work has been partially funded by the EU H2020 \textit{Energy oriented Center of Excellence (EoCoE) for computer applications}, the \textit{New European Wind Atlas (NEWA)}, the \textit{High Performance Computing for Energy (HPC4E)}, and the \textit{Supercomputing and Energy for Mexico (ENERXICO)} projects. We thank Iberdrola Renovables \cite{IberdrolaRenovables} for their collaboration and for providing the topography data of the different scenarios.


\end{document}